\documentclass{amsart}

\usepackage{latexsym,amssymb,amsmath,amscd,amsthm,amsxtra}
\usepackage{amsmath}
\usepackage{amsfonts}
\usepackage{amssymb}
\usepackage{graphicx}

\newtheorem{thm}{Theorem}[section]
\newtheorem{prop}[thm]{Proposition}
\newtheorem{lem}[thm]{Lemma}
\newtheorem{cor}[thm]{Corollary}
\newtheorem{dfn}[thm]{Definition}

\newtheorem{rmk}[thm]{Remark}
\newtheorem{ex}[thm]{Example}

\newtheorem{conj}[thm]{Conjecture}

\newcommand{\complex}{{\mathbb C}}

\newcommand{\reals}{{\mathbb R}}

\newcommand{\integers}{{\mathbb Z}}

\newcommand{\Hom}{{\rm Hom}}
\newcommand{\Poly}{{\rm Poly}}

\newcommand{\calm}{{\mathcal M}}

\begin{document}

\title{Noncommutative Poisson Structures on Orbifolds}

\author{Gilles Halbout}
\address{Institut de Math\'ematiques et de Mod'elisation de Montpellier I3M, UMR 5149, Universit\'e de Montpellier 2, F-34095 Montpellier cedex 5, France}
\email{ghalbout@darboux.math.univ-montp2.fr}

\author{Xiang Tang}
\address{Department of Mathematics, Washington University, St.
   Louis, MO, 63130, USA}
\email{xtang@math.wustl.edu}
\thanks{The second author was supported in part by NSF Grant 0604552.}


\date{}



\begin{abstract}
In this paper, we compute the Gerstenhaber bracket on the Hoch-schild
cohomology of $C^\infty(M)\rtimes G$ for a finite group $G$ acting
on a compact manifold $M$. Using this computation, we obtain
geometric descriptions for all  noncommutative Poisson structures on
$C^\infty(M)\rtimes G$ when $M$ is a symplectic manifold. We also
discuss examples of deformation quantizations of these
noncommutative Poisson structures.
\end{abstract}

\maketitle

\section{Introduction}
It is well known \cite {ge:def-hoc} that the deformation theory of
an associative algebra $A$ is closely related to the Hochschild
cohomology $HH^\bullet (A;A)$ of $A$. In particular, the
infinitesimal deformation of $A$ is governed by $HH^2(A;A)$.
Furthermore, if we want the infinitesimal deformation to be
integrable, it is necessary (but not sufficient) to
require that the 2-cocycle $\Pi\in C^2(A;A)$ associated to the
infinitesimal deformation satisfies the equation $[\Pi, \Pi]_G=0$
in $HH^3(A;A)$, where $[\ ,\ ]_G$ is the Gerstenhaber bracket on
$HH^\bullet (A;A)$.

Suppose that $A$ is an algebra of smooth functions on a smooth manifold
$M$. Then, according to the Hochschild-Kostant-Rosenberg theorem,
the second Hoch-schild cohomology classes in $HH^2 (A;A)$ satisfying the
above integrability conditions are in one to one correspondence with
Poisson structures on $M$. Inspired by this relationship between
Poisson geometry and deformation theory, Jonathan Block and Ezra
Getzler \cite{bl-ge:quantization} and Ping Xu \cite{xu}
independently introduced a notion of a noncommutative Poisson
structure on an associative algebra in early 90's.

\begin{dfn}
\label{dfn:nc-Poisson} A noncommutative Poisson structure on an
associative algebra $A$ is an element $\Pi$ in the second
Hochschild cohomology group $H^2(A,A)$ of $A$, whose Gerstenhaber
bracket with itself vanishes, i.e. $[\Pi, \Pi]_G=0$.
\end{dfn}

In this paper, we want to study noncommutative Poisson structures on
orbifolds coming from global quotients. Let $M$ be a compact smooth
manifold, and $G$ be a finite group acting on $M$. (For simplicity, we will always assume in this paper that the $G$-action on $M$ is effective.)  Our orbifold is
the quotient space $X=M/G$. Because $X$ is usually a topological
space with quotient singularities, the algebra $C^\infty(M)^G$ of
$G$-invariant smooth functions on $M$ is not regular. Taking a lesson
from noncommutative geometry \cite{c:koszul}, we consider the
crossed product algebra $C^\infty(M)\rtimes G$ as a natural
replacement. The crossed product algebra $C^\infty(M)\rtimes G$, thought it is noncommutative, has very nice algebraic properties.  Our main aim in this paper is to find out geometric
descriptions for all noncommutative Poisson structures on
$C^\infty(M)\rtimes G$ when $M$ is a symplectic manifold with a symplectic action.

As a first step toward understanding noncommutative Poisson
structures on $C^\infty(M)\rtimes G$, the second author and
the his coauthors \cite{nppt} computed the Hochschild cohomology of
$C^\infty(M)\rtimes G$ as a vector space:
\begin{equation}
\label{eq:cohomology} HH^\bullet (C^\infty(M)\rtimes G;
C^\infty(M)\rtimes G)=\Gamma^\infty(\bigoplus_{g\in
G} \wedge^{\bullet -l(g)}TM^g\otimes
\wedge^{l(g)}N^g)^G.
\end{equation}
We explain the notation in the above equation. First, $M^g$ is the fixed point manifold of $g$, and $N^g$ is the normal bundle of the embedding of $M^g$ in $M$, where $l(g)$ is  the dimension of $N^g$. The group $G$ acts on disjoint
union $\sqcup M^g$ of $M^g$ for all $g\in G$ by the conjugate action. We
remark that $M^g$ may have different components with different
dimensions, and we take the disjoint union of all the components,
and accordingly $l(g)$ is taken as a local constant function on
$M^g$. (Following \cite{nppt}, in this paper we view
$C^\infty(M)\rtimes G$ as a bornological algebra with the bornology
defined by the Frech\'et topology. We take $HH^\bullet$ to be the
continuous Hochschild cohomology of a bornological algebra.
Accordingly, all computations and constructions in this paper are
local respect to the orbifold $M/G$. We often work with a vector
space (or a $G$-invariant open subset) with a linear $G$ action,
which we refer to as ``local'' computation.) We call the stratified
space $\sqcup M^g/G$ the inertia orbifold $\tilde{X}$ associated to
$X=M/G$. In other words, we can interpret Equation (\ref{eq:cohomology})
 as saying that the Hochschild cohomology of $C^\infty(M)\rtimes G$ is equal to the
space of ``vector fields" on $\tilde{X}$.

The main difficulty in finding noncommutative Poisson structures
on $C^\infty(M)\rtimes G$ is to compute the Gerstenhaber
bracket on $HH^\bullet (C^\infty(M)\rtimes G;
C^\infty(M)\rtimes G)$.
To compute the Gerstenhaber bracket, we need
to have quasi-isomorphisms between the Hochschild cochain complex
\[
C^\bullet (C^\infty(M)\rtimes G; C^\infty(M)\rtimes G)
\]
and
\[
\Gamma^\infty(\bigoplus_{g\in G} \wedge^{\bullet
-l(g)}TM^g\otimes \wedge^{l(g)}N^g)^G.
\]
In \cite{nppt}, we defined a quasi-isomorphism $L$ in the
following direction:
\[
L: C^\bullet (C^\infty(M)\rtimes G; C^\infty(M)\rtimes
G)\longrightarrow \Gamma^\infty(\bigoplus_{g\in G}
\wedge^{\bullet -l(g)}TM^g\otimes
\wedge^{l(g)}N^g)^G.
\]
The first step of this paper is to define a quasi-isomorphism $T$
in the other direction. This turns out to be much harder to
construct than the map $L$ in \cite{nppt}. We need to construct some nonlocal operators on
$C^\infty(M)$, which we call twisted cocycles. These cocycles are
closely related to the Lusztig-Demazure operator ({\it cf.}
\cite{ob:lus-dem}). Using the maps $T$ and $L$, we are able to
compute the Gerstenhaber brackets on $HH^\bullet (C^\infty(M)\rtimes
G; C^\infty(M)\rtimes G)$. Our results show that the Gerstenhaber
bracket on orbifolds is a generalization of the classical
Schouten-Nijenhuis bracket on manifolds. We call this bracket the
twisted Schouten-Nijenhuis bracket on
$\Gamma^\infty(\oplus_{g}\wedge^{\bullet-l(g)} TM^g \otimes
\wedge^{l(g)}N^g)^G$. Using the twisted Schouten-Nijenhuis bracket,
we are able to solve the equation $[\Pi, \Pi]_G=0$ on
$HH^2(C^\infty(M)\rtimes G; C^\infty(M)\rtimes G)$ geometrically.
This leads to a full description of noncommutative Poisson
structures on $C^\infty(M)\rtimes G$ when $M$ is a symplectic
manifold.

If we consider a complex symplectic vector space $V$ with a
symplectic $G$ action, the cocycles used in the definition of
symplectic reflection algebras \cite{eg:cherednik} correspond to a
special class of noncommutative Poisson structures on
$\Poly(V)\rtimes G$, where $\Poly(V)$ is the algebra of polynomials
on $V$. Using the results from \cite{eg:cherednik} we prove in this
paper that all these cocycles can be extended to a formal
deformation of the algebra $\Poly(V)\rtimes G$. As a generalization,
we expect that all the noncommutative Poisson structures discovered
in this paper can be extended to formal deformations, which will
generalize the symplectic reflection algebras. This question is
closely related to the following formality conjecture on orbifolds,
which we will discuss in future publications.
\begin{conj}
The Hochschild complex of the algebra $C^\infty(M)\rtimes G$ is a
formal differential graded Lie algebra.
\end{conj}
In the last part of this paper, we provide concrete new examples
of noncommutative Poisson structures on $Poly(\reals^4)\rtimes
(\integers_n\times \integers_m)$ with
$\integers_n=\integers/\integers_n$ and
$\integers_m=\integers/m\integers$. These Poisson structures are
not symplectic at all, and instead should be viewed as
noncommutative quadratic Poisson structures. The connection
between these ``noncommutative quadratic Poisson structures" and
quantum R matrices will be studied in the near future. In general,
there are many interesting examples of noncommutative Poisson
structures on orbifolds. We are working with Jean-Michel Oudom in \cite{hot:example} on this material.

Besides the Gerstenhaber bracket, there is also a product on the
Hochschild cohomology $HH^\bullet (C^\infty(M)\rtimes G;
C^\infty(M)\rtimes G)$. In \cite{pptt:product}, with Pflaum,
Posthuma and Tseng, the second author studied the ring structure on
the Hochschild cohomology of the deformed algebras of
$C^\infty(M)\rtimes G$, which is closely related to the Chen-Ruan
orbifold cohomology \cite{cr:orbifold}.

This paper is organized as follows. In Section 2, we briefly recall
the basic definitions in Hochschild cohomology and the computations
of $HH^\bullet(C^\infty(M)\rtimes G, C^\infty(M)\rtimes G)$ as a
vector space in \cite{nppt}. In Section 3, we will focus on the
construction of twisted cocycles and a quasi-isomorphism $T$,
\[
T: \Gamma^\infty(\bigoplus_{g\in G} \wedge^{\bullet
-l(g)}TM^g\otimes
\wedge^{l(g)}N^g)^g\longrightarrow C^\bullet
(C^\infty(M)\rtimes G; C^\infty(M)\rtimes G).
\]
In Section 4, we will study the Gerstenhaber bracket on the
Hochschild cohomology $HH^\bullet (C^\infty(M)\rtimes G;
C^\infty(M)\rtimes G)$. In Section 5, we will give a full
description of noncommutative Poisson structures on
$C^\infty(M)\rtimes G$ when $M$ is a symplectic manifold.
Furthermore we will study formal deformation of a special type of
noncommutative Poisson structures and compute their second Poisson
cohomology group. We end this section by showing two new families of
noncommutative quadratic Poisson structures on
$\complex^2/\integers_n\times \integers_m$.

\begin{rmk}
Unless otherwise specified, we work with the field $\reals$, real
vector spaces, and real manifolds. However, many results in this paper have
analogs with the field $\complex$, complex vector space and affine
varieties.
\end{rmk}

\noindent{\bf Acknowledgments}: We would like to thank Vassiliy
Dolgushev, Benjamin Enriquez, Pavel Etingof, Victor Ginzburg,
Jean-Michel Oudom, Hessel Posthuma, and Markus Pflaum for useful
conversations. We would like also to particularly thank Ping Xu for
hosting our visits to Penn State University in Spring 2005, where we
started this project.

\section{Hochschild cohomology of $C^\infty(M)\rtimes G$}
In this section, we briefly review our work in \cite{nppt} on computing the Hochschild cohomology of the algebra $C^\infty(M)\rtimes
G$ as a vector space.

\subsection{Hochschild cohomology of an algebra}
We review in this subsection the definition of Hochschild cohomology of an algebra. Let
$A$ be a unital algebra over $\reals$, and $M$ be a bimodule of $A$. The degree $n$
Hochschild cochain complex $C^n(A, M)$ of $A$ with coefficients in
$M$ consists of $\reals$-linear maps from $A^{\otimes n}$ to $M$,
i.e. $\Hom(A^{\otimes n},M)$. The coboundary map
$\partial$ on Hochschild cochains $C^\bullet(A, M)$ is defined by
\begin{eqnarray*}
\partial:&&C^n(A, M)\longrightarrow C^{n+1}(A,M)\\
\partial(\xi )(a_0, \cdots, a_n)&=& a_0\xi(a_1, \cdots, a_n)+\sum_{i=1}^n(-1)^i\xi (\cdots, a_{i-1}a_i, \cdots )\\
&&\qquad \qquad +(-1)^{n+1}\xi(a_1, \cdots, a_{n-1})a_n.
\end{eqnarray*}
It is a straightforward check that $\partial^2=0$.

\begin{dfn}\label{dfn:hoc}
Define the Hochschild cohomology of $A$ with coefficient in $M$ to
be the cohomology of the differential cochain complex $(C^\bullet(A,
M), \partial)$.
\end{dfn}

\begin{rmk}
\label{rmk:hoch}
In the case that $A$ is
a topological (or bornological) algebra, then we need to consider a
topological (bornological) algebra bimodule $M$, and the topological
(bornological) tensor product of $A$, and continuous (bounded)
linear maps from $A$ to $M$. All the definitions we introduce below
naturally extend to topological (bornological) version and we refer
to \cite{pptt:product}[Appendix A] for details. The above definitions work for arbitrary fields.
\end{rmk}

The left and right multiplication of $A$ on $A$ makes $A$ a bimodule
of $A$. Gerstenhaber \cite{ge:def-hoc} used the Hochschild
cohomology $H^\bullet(A,A)$ to study deformation theory of $A$. On
$C^\bullet(A,A)$, besides the coboundary differential, there are two
interesting operations, the
\begin{enumerate}
\item Gerstenhaber bracket, and the
\item cup product.
\end{enumerate}

Then Gerstenhaber bracket is used in defining a noncommutative Poisson
structure in Definition \ref{dfn:nc-Poisson}. We recall its
definition.
\begin{dfn}\label{dfn:ger}
We define a pre-Lie
product $\circ$ on $C^\bullet (A;A)$. For $\phi\in C^k(A;A)$,
$\psi\in C^l(A;A)$, $\phi\circ \psi\in C^{k+l-1}(A;A)$ is defined
by
\[
\begin{split}
&\phi\circ\psi(a_1, \dots ,
a_{k+l-1})\\
=&\sum_{i=1}^{k}(-1)^{(i-1)(l-1)}\phi(a_1, \dots ,
a_{i-1}, \psi(a_i, \dots , a_{i+l-1}), a_{i+l}, \dots ,
a_{k+l-1}),
\end{split}
\]
for $a_i\in A$, $i=1, \dots  k+l-1$. The Gerstenhaber bracket $[\
,\ ]$ on $C^\bullet (A;A)$ is defined to be the commutator of the
pre-Lie product, i.e.,
\[
[\phi, \psi]_G=\phi\circ \psi-(-1)^{(k-1)(l-1)}\psi\circ \phi.
\]
\end{dfn}
The Gerstenhaber bracket is compatible with the differential on
$C^\bullet(A,A)$, and therefore defines a Lie bracket on the
Hochschild cohomology $HH^\bullet(A,A)$.
\subsection{$HH^\bullet(Poly(V)\rtimes G, Poly(V)\rtimes G)$}

In this subsection, we consider the algebra $C^\infty(M)\rtimes G$
for a finite group action on a compact manifold $M$, and explain the
computation in \cite{nppt} of the Hochschild cohomology
$HH^\bullet(C^\infty(M)\rtimes G, C^\infty(M)\rtimes G)$. In this
subsection we will mainly focus on the local case, i.e., $G$-linear action
on a vector space $V$; in the next subsection, we explain
the extension of the computations general manifolds.

Let $Poly(V)$ be the algebra of polynomial functions on a vector
space $V$, and $G$ be a finite group acting linearly on $V$. We
computed \cite{nppt} the Hochschild cohomology of the crossed
product algebra $Poly(V)\rtimes G$.  The major step is a
construction of a quasi-isomorphism
\[
L: C^\bullet (Poly(V)\rtimes G, Poly(V)\rtimes G)\longrightarrow
\Gamma^{\infty}(\bigoplus_{g\in G}\wedge^{\bullet
-l(g)}TV^g\otimes \wedge^{l(g)}N^g)^{G}.
\]
This map was constructed implicitly in the proof of Theorem 3.1,
\cite{nppt}. We make it explicit in the following.

The map $L$ is a composition of three cochain maps:

\noindent{\bf Step I}:
\[
L_1: C^\bullet (Poly(V)\rtimes G, Poly(V)\rtimes G)\longrightarrow
\big(C^\bullet (Poly(V), Poly(V)\rtimes G)\big)^G,
\]
where $G$ acts on $C^\bullet (Poly(V), Poly(V)\rtimes G)$ by
\[
g\Psi(a_1, \dots, a_n)=U_{g^{-1}}\cdot \Psi(g(a_1),
\dots, g(a_n))\cdot U_{g}.
\]

Here $U_{g}$ denotes the element $g$ seen in $Poly(V)\rtimes G$.
Given a Hochschild cocycle $\Psi\in C^k(Poly(V)\rtimes G,
Poly(V)\rtimes G)$, we define $L_1(\Psi)\in \big(C^k(Poly(V), Poly(V)\rtimes
G)\big)^G$,
\[
L_1(\Psi)(f_1, \dots , f_k)=\frac{1}{|G|}\sum_g (g\Psi)(f_1, \dots ,
f_k), \ \ \ \forall f_1, \dots , f_k\in Poly(V),
\]
where $|G|$ is the order of group $G$.

\noindent{\bf Step II}:
 \[
L_2: \big(C^\bullet (Poly(V), Poly(V)\rtimes
G)\big)^G\longrightarrow \big(\bigoplus_{g\in G}
\Gamma^\infty(\wedge^\bullet TV), \kappa_g\wedge\ \big)^G.
\]

Let $A_g$ be a vector space isomorphic to $Poly(V)$ equipped with the $Poly(V)$-bimodule structure
\[
a\cdot \xi \cdot (b)=a\xi g(b),\ \text{for}\ a,b\in Poly(V),\ \xi\in Poly(V)_g,
\]
where the right hand side is the product of $a, \xi$, and $g(b)$ as elements
in $Poly(V)$. As a $Poly(V)$-$Poly(V)$ bimodule, $Poly(V)\rtimes G$
has a natural splitting into a direct sum of submodules
$\oplus_{g\in G}A_g$. Correspondingly, the cochain complex
$C^\bullet (Poly(V), Poly(V)\rtimes G)$ has a natural splitting into
$\oplus_{g\in G}C^\bullet (Poly(V), A_g)$. We define $L_2$ to be the
sum of the maps
\[
L_2^g: C^\bullet (Poly(V), A_g)\longrightarrow
\big(\Gamma^\infty(\wedge^\bullet TV), \kappa_g \big)
\]
over all $g\in G$.

On $V$, we introduce the vector field
$X(x)=\sum_{i}x^i\frac{\partial}{\partial x^i}$, where the $x^i$ are
coordinate functions on $V$. Define the vector field $\kappa_g\in
\Gamma^\infty(TV)$ by
\[
\kappa_g(x)=X(g(x))-X(x).
\]

We notice that for a permutation $\sigma$ of $k$ elements fixing
$x\in V$, the product
$(x_{\sigma(1)}^{i_1}-x^{i_1})\cdots(x_{\sigma(k)}^{i_k}-x^{i_k})$
is a function on $x_1, \cdots, x_k\in V$ by taking the product of
the values of the coordinate functions. Given an element $\Psi\in
C^k(Poly(V), A_g)$, we define $L_2^g(\Psi)\in
\Gamma^\infty(\wedge^kTV)$, the usual projection to anti-symmetric
linear operators, by
\[
L_2^g(\Psi)(x)=\sum\limits_{i_1, \dots , i_k}\frac{1}{k!}\sum_{\sigma\in S_k}(-1)^\sigma
\Psi\big((x_{\sigma(1)})^{i_1}-x^{i_1}) \cdots
(x_{\sigma(k)}^{i_k}-x^{i_k})\big)(x)\frac{\partial}{\partial x^{i_1}}\wedge
\cdots \wedge \frac{\partial}{\partial x^{i_k}},
\]
where $S_k$ is the permutation group of $k$-elements.

The $G$ action on $\big(\bigoplus_{g\in G} \Gamma^\infty(\wedge^\bullet TV),
\kappa_g\wedge\ \big)$ is defined by
\[
h\big(\sum_g\phi_g\big)|_{h^{-1}gh}=h_*(\phi_g),\qquad \text{for }\sum_g\phi_g \in \bigoplus_g\Gamma^\infty(\wedge^\bullet TV), \ \text{and\ }h\in G.
\]

It is straightforward to check that $L_2^g$ is $G$-equivariant, and
therefore defines a map
\[
L_2:\big(C^\bullet (Poly(V), Poly(V)\rtimes G)\big)^G\longrightarrow
\big(\bigoplus_{g\in G} \Gamma^\infty(\wedge^\bullet TV),
\kappa_g\wedge\ \big)^G.
\]

\noindent{\bf Step III}:
\[
L_3:\Big(\bigoplus_{g\in G} \big( \Gamma^\infty(\wedge^\bullet TV),
\kappa_g\wedge\ \big)\Big)^G\longrightarrow \Big(\bigoplus_{g\in
G}\big(\Gamma^\infty(\wedge^{\bullet -l(g)} TV^g\otimes
\wedge^{l(g)}N^g), 0\ \big)\Big)^G.
\]

Let $C_g$ be the cyclic group generated by $g$, which has a natural action on $V$. As $C_g$ is abelian, $V$ is decomposed into a direct sum of $C_g$ irreducible representations. Let $V^g$ be the subspace of  all trivial
$C_g$-representations in $V$, and $N^g$ be the sum of all nontrivial irreducible $C_g$ representations in $V$. Therefore, $V$ can be written as $V^g\oplus N^g$.

We define $L_3$ to be the sum of $L_3^g$, which is defined to
be
\[
L_3^g(X)=pr^g( X|_{V^g}),
\]
where $X|_{V^g}$ is the restriction of $X\in \wedge^\bullet TV$ to
$\wedge^\bullet TV|_{V^g}$, and $pr^g$ projects $\wedge^\bullet
TV|_{V^g}$ to $\wedge^{\bullet -l(g)} TV^g\otimes \wedge^{l(g)}N^g$.

The space $\bigoplus_{g\in
G}\big(\Gamma^\infty(\wedge^{\bullet -l(g)} TV^g\otimes
\wedge^{l(g)}N^g, 0 \big)$ is closed under $G$ action on $\big(\bigoplus_{g\in G} \Gamma^\infty(\wedge^\bullet TV),
\kappa_g\wedge\ \big)$ and therefore inherits a $G$ action. Similarly to the computations taken in Step II, we can easily check that $L_3^g$ is
$G$-equivariant and that therefore $L_3$ defines a map on the $G$-invariant
components.

We proved in Section 3 of \cite{nppt} that $L=L_3\circ L_2\circ L_1$
is a quasi-isomorphism of cochain complexes
\[
L: C^\bullet (Poly(V)\rtimes G, Poly(V)\rtimes G)\longrightarrow
\Gamma^{\infty}\big(\bigoplus_{g\in G}\wedge^{\bullet
-l(g)}TV^g\otimes \wedge^{l(g)}N^g\big)^{G}.
\]
\subsection{The general cases}
In the above  Steps I-III, we have explained how to compute the
Hochschild cohomology of $Poly(V)\rtimes G$ using the
quasi-isomorphism $L$. We explain briefly how to generalize this
construction to $C^\infty(M)\rtimes G$ by defining $L$ to be
\[
L:C^\bullet(C^\infty(M)\rtimes G, C^\infty(M)\rtimes
G)\longrightarrow \Gamma^\infty\big(\bigoplus_{g\in
G}\wedge^{\bullet-l(g)}TM^g\otimes \wedge^{l(g)}N^g\big)^G.
\]

Firstly, we observe that the map $L_1$ is a quasi-isomorphism
from $C^\bullet(A\rtimes G, A\rtimes G)$ to $\big(C^\bullet(A,
A\rtimes G)\big)^G$ that is true for any algebra $A$ with a finite group
action. Therefore, the map $L_1$ extends to the general case
$C^\infty(M)\rtimes G$ naturally:
\[
L_1: C^\bullet(C^\infty(M)\rtimes G, C^\infty(M)\rtimes
G)\longrightarrow \Big(C^\bullet\big(C^\infty(M), C^\infty(M)\rtimes
G\big)\Big)^G.
\]

Secondly, $L_3$ is generalized to the manifold case, as the map
\[
L_3:\Big(\bigoplus _{g\in G}\big(\Gamma^\infty(\wedge^\bullet TM),
\kappa \wedge \big)\Big)^G\longrightarrow \Big(\bigoplus_{g\in
G}\big(\Gamma^\infty(\wedge^{\bullet-l(g)}TM^g\otimes
\wedge^{l(g)}N^g), 0\big)\Big)^G.
\]
obtained by composing the projection map $\wedge^\bullet TM|_{M^g}\to
\wedge^{\bullet-l(g)}TM^g\otimes \wedge^{l(g)}N^g$ with the
restriction map from $\wedge^\bullet TM$ to $\wedge^\bullet TM|_V$.

Thirdly, to generalize $L_2$, we use Connes' map
\cite[Lemma 44]{c:koszul} from the Koszul resolution of $C^\infty(M)$ to
its Bar resolution for a manifold with an affine structure. We can then define $L_2^g$ from the Hochschild cochain complex $C^\bullet
(C^\infty(M), C^\infty(M)_g)$ to $\wedge^\bullet  TM$ in a similar
way. There is a small issue here that the Connes construction
\cite[Lemma 44]{c:koszul} only works for affine manifolds.
Therefore, for a general manifold $M$, we need to realize the
Hochschild cochain complex of $C^\infty(M)\rtimes G$ as a (pre)sheaf
over the orbifold $M/G$, and use \v Cech techniques to compute the
sheaf cohomology of this (pre)sheaf. In this framework, $L$ will be
a quasi-isomorphism of (pre)sheaves which is locally defined
as the $L_g$ was in Section 2.2, Step II. We refer to
\cite[Section 3]{pptt:product} for details.

At the end of this subsection, we explain the following observation
which makes our study easier. Recall Equation (\ref{eq:cohomology}) concerning the Hochschild cohomology,
\[
HH^\bullet (C^\infty(M)\rtimes G; C^\infty(M)\rtimes
G)=\Gamma^\infty(\bigoplus_{g\in G} \wedge^{\bullet
-l(g)}TM^g\otimes \wedge^{l(g)}N^g)^G.
\]
As $g$ is in the $g$-centralizer subgroup of $G$, the $g$-fixed
point component's contribution to $HH^\bullet( C^\infty(M)\rtimes G,
C^\infty(M)\rtimes G)$ has to be from $g$-invariant sections of
$\Gamma^\infty(\wedge^{\bullet-l(g)}TV^g\otimes \wedge^{l(g)}N^g)$.
As $g$ acts on $TV^g$ trivially, a $g$-invariant section of
$\wedge^{\bullet-l(g)}TV^g\otimes \wedge^{l(g)}N^g$ must have
$g$-invariant component in $\wedge^{l(g)}N^g$. Note that $\wedge^{l(g)}N^g$ is
a line bundle over $V^g$. To have a nonzero $g$-invariant section,
we need that the $g$ action on $\wedge^{l(g)}N^g$ must be trivial. This
implies that $\det(g|_{N^g})=1$. Furthermore, we notice that
$g|_{N^g}$'s action on $N^g$ is of finite order and therefore can be
diagonalized. If $N^g$ is of odd dimension, then by the fact that
$\det(g|_{N^g})=1$, we conclude that $g|_{N^g}$ has at least one
eigenvalue equal to 1. This contradicts to the assumption of $N^g$.
Therefore, we conclude that if $\dim(N^g)$ is odd, there is no
nonzero contribution to $HH^\bullet(C^\infty(M)\rtimes G,
C^\infty(M)\rtimes G)$ from this $g$-fixed point component. Hence,
the Hochschild cohomology of $C^\infty(M)\rtimes G$ has no
contribution from $g$-fixed point submanifolds with odd $l(g)$.
Therefore, we conclude
\begin{equation}
\label{eq:coh-even}HH^\bullet (C^\infty(M)\rtimes G;
C^\infty(M)\rtimes G)=\Gamma^\infty\Big(\bigoplus_{\small
\begin{array}{c}g\in
G,\\
 l(g) \text{ is even}\end{array}} \wedge^{\bullet -l(g)}TM^g\otimes \wedge^{l(g)}N^g\Big)^G.
\end{equation}

\section{Hochschild cohomology and Quasi-isomorphisms}
In this section, we construct a quasi-isomorphism
\[
T: \Gamma^\infty\big(\bigoplus_{g\in G} \wedge^{\bullet
-l(g)}TM^g\otimes \wedge^{l(g)}N^g\big)^G\longrightarrow C^\bullet
(C^\infty(M)\rtimes G; C^\infty(M)\rtimes G),
\]
which is a quasi-inverse to the map $L$ reviewed in Sections 2.2-2.3.
First, we construct a twisted cocycle $\Omega_g$ for each
element $g$ associated to the determinant line bundle
$\wedge^{l(g)}N^g$; second, we use the twisted cocycles to
construct the map $T$. We mainly focus on the local case
$Poly(V)\rtimes G$, and explain at the end how to generalize the
construction to general manifolds.
\subsection{Twisted cocycle}
To construct the map $T$, we need to understand cocycles in
$C^\bullet(Poly(V),  A_g)$, which we call $g$-twisted cocycles.
(Recall that $A_g$ is a bimodule of $Poly(V)$ isomorphic to
$Poly(V)$ as a vector space but with the $g$-twisted multiplication
from the right (Section 2.2, Step II).) If we trace through the
computation Steps I-III in Section 2.2, we see that
$H^\bullet(Poly(V), A_g)$ is isomorphic to
$\Gamma^\infty(\wedge^{\bullet-l(g)}TV^g\otimes \otimes
\wedge^{l(g)}N^g)$. To construct the quasi-inverse of $L$, we need to
associate to each element in
$\Gamma^\infty(\wedge^{\bullet-l(g)}TV^g \otimes \wedge^{l(g)}N^g)$
a cocycle in $C^\bullet(Poly(V),  A_g)$. In particular, we need
to associate to  an element in $\wedge^{l(g)}N^g$ a degree $l(g)$
cocycle in $C^\bullet(Poly(V),  A_g)$. It is well-known that an
element in $\Gamma^\infty(\wedge^{\bullet-l(g)}TV^g\otimes
\wedge^{l(g)}N^g)$ can be viewed as a multi-differential operator on
$V$. A natural guess is that this multi-differential operator might
be a cocycle in $C^\bullet(Poly(V),  A_g)$. However, one can easily
check that such a multi-differential operator is not a cocycle in
$C^\bullet(Poly(V), A_g)$ except in the case $g=1$. (For instance, because the
right module structure on $A_g$ is twisted by the $g$ action, a
degree 1 cocycle $\xi$ in $C^\bullet(Poly(V), A_g)$ must satisfy
$\xi(ab)=a\xi(b)+\xi(a)g(b)$. One can quickly check that when $g\ne 1$, no nonzero
vector fields on $V$ can satisfy such an equation.)
Therefore, we have to modify the above natural guess so that the
outcome is a cocycle in $C^\bullet(Poly(V),  A_g)$. This leads us to
the following construction.

Recall that $V=V^g\oplus N^g$. We choose coordinates $x^1, \cdots,
x^{n-l(g)}, y^1,\cdots, y^{l(g)}$ on $V$ ($n=\dim(V)$), such that
$x^{1}, \dots , x^{n-l(g)}$ are coordinates on $ V^g$ and $y^1,
\dots , y^{l(g)}$ are coordinates on $N^g$. And we assume that $g$
action on $y^{1}, \dots, y^{l(g)}$ is diagonalized\footnote{We can
pass to the complex coordinates to achieve such a choice. }.

For any $\sigma \in S_{l(g)}$, the permutation group of $l(g)$
elements, we introduce the following vectors in $N^g$. Let
$(\tilde{y}^1, \dots , \tilde{y}^{l(g)})=g(y^{1}, \dots$ $
y^{l(g)})$. Define
\[
\begin{split}
z^0_\sigma=(y^1, \dots , y^{l(g)})\ \ \ \ \ \ \ \ \ \ \ &
z^1_\sigma=(y^1, \dots , \tilde{y}^{\sigma(1)}, \dots )\\
z^2_\sigma=(y^1, \dots , \tilde{y}^{\sigma(1)}, \dots ,
\tilde{y}^{\sigma(2)}, \dots ,)\ \ \ \ &
\cdots\\
z^{l(g)-1}_\sigma=(\tilde{y}^1, \dots ,
y^{\sigma(l(g))}, \dots ,)\ \ \ \ \ \
&z^{l(g)}_\sigma=(\tilde{y}^1, \dots ,
\tilde{y}^{l(g)}).
\end{split}
\]

Let $\Omega_{g}$ be a  $l(g)-$cochain in $C^\bullet (A,
A_g)$ defined as follows: {\small
\begin{equation}
\begin{array}{cc}
\label{eq:twisted-cocycle} &\Omega_{g}(f_1, \dots ,
f_{l(g)})(x,
y)\\
&:=\frac{1}{l(g)!}\sum\limits_{\sigma\in
S_{l(g)}}(-1)^\sigma\frac{(f_1(x,z_\sigma^0)-f_1(x,
z_\sigma^1))(f_2(x, z_\sigma^1)-f_2(x, z_\sigma^2))\cdots
(f_{l(g)}(x,z_\sigma^{l(g)-1})-f_{l(g)}(x,
z_\sigma^{l(g)}))}{(y^1-\tilde{y}^1)\cdots
(y^{l(g)}-\tilde{y}^{l(g)})}.
\end{array}
\end{equation}
}

We remark that when $\sigma$ is the identity permutation, we
have
\[
\begin{split}
z^0_{id}=(y^1, \dots , y^{l(g)})\ \ \ \ \ \ \ \ \ \ \ &
z^1_{id}=(\tilde{y}^1, \dots , y^{l(g)} )\\
\cdots& \cdots\\
z^{l(g)-1}_{id}=(\tilde{y}^1, \dots ,
\tilde{y}^{l(g)-1}, y^{l(g)})\ \ \ \ \ \
&z^{l(g)}_{id}=(\tilde{y}^1, \dots , \tilde{y}^{l(g)}).
\end{split}
\]

The corresponding contribution in the summation of expression
(\ref{eq:twisted-cocycle}) is {\tiny
\[
\frac{(f_1(x,z_{id}^0)-f_1(x, z_{id}^1))(f_2(x, z_{id}^1)-f_2(x,
z_{id}^2))\cdots
(f_{l(g)}(x,z_{id}^{l(g)-1})-f_{l(g)}(x,
z_{id}^{l(g)}))}{(y^1-\tilde{y}^1)\cdots
(y^{l(g)}-\tilde{y}^{l(g)})},
\]}
which converges to
\[
\frac{\partial}{\partial y^1}f_1(x,0)\cdots
\frac{\partial}{\partial y^{l(g)}}f_{l(g)}(x,0),
\]
as $y^1,\cdots, y^{l(g)}$ go to 0.

Therefore the identity component in Equation
(\ref{eq:twisted-cocycle}) can be viewed as a $g$ analog of
the multi-differential operator
\[
\frac{\partial}{\partial y^1}\otimes \cdots \otimes
\frac{\partial}{\partial y^{l(g)}}.
\]

Summing over all permutations, $\Omega_g$ can be viewed as a
$g$-analog of the multi-differential operator
\[
\Lambda_g=\frac{\partial}{\partial y^1}\wedge\cdots \wedge
\frac{\partial}{\partial y^{l(g)}}.
\]

It is straightforward to check that $\Omega_g$ is a cocycle in
$C^{l(g)}(A, A_g)$.

The following is a simple property of $\Omega_g$.
\begin{lem}
\label{lem:linear-twisted}The twisted cocycle $\Omega_g$ satisfies
the following properties:
\begin{enumerate}
\item
\[
\Omega_g(y^1\otimes\cdots\otimes y^{l(g)})=1;
\]
\item
\[
\Omega_g(y^1\otimes\cdots \otimes c\otimes \cdots \otimes
y^{l(g)})=0,
\]
when $c$ is a constant function.
\end{enumerate}
\end{lem}
\begin{proof}
A straightforward check.
\end{proof}

We remark that the expression (\ref{eq:twisted-cocycle}) of
$\Omega_g$ depends on the choices of coordinates $y^1, \cdots,
y^{l(g)}$. This makes $\Omega_g$ in general not invariant with
respect to the conjugate action. However, we have the following
property.

\begin{prop}
\label{prop:conjugation-twist-cocy}Let $C(g)$ be the centralizer
subgroup of $g$, which acts on $N^g$. If $C(g)$ action on $N^g$ is
diagonalizable\footnote{By a diagonalizable action, we mean $C(g)$
action on $N^g$ splits into a direct sum of 1-dim or 2-dim
representations of $C(g)$.}, there is a natural construction of
$\Omega_g$ such that
\[
h(\Omega_g)=\det(h|_{N^g})\Omega_g, \qquad h\in C(g).
\]
\end{prop}

\begin{proof}
As $C(g)$ action on $N^g$ is diagonalizable and $g$ commutes with
elements in $C(g)$, $g$ and $C(g)$ action on $N^g$ can be
diagonalized simultaneously. Therefore, we can find coordinates
$y^1, \cdots, y^{l(g)}$ on $N^g$, which are eigenfunctions of $g$
and $C(g)$ action. We define $\Omega_g$ using the coordinates $y^i$
as Equation(\ref{eq:twisted-cocycle}). In particular,
$\tilde{y}^i=g^iy^i$, and $h(\tilde{y}^i)=h(g^iy^i)=g^ih^iy^i$,
where $g^i$ and $h^i$ are eigenvalues of $g$ and $h$ action on
$y^i$. Plug these expressions in the definition of $h(\Omega_g)$, we
obtain the equation
\[
h(\Omega_g)=\det(h|_{N^g})\Omega_g, \qquad h\in C(g).
\]
\end{proof}

There are two special cases we know that the conditions assumed in
Proposition \ref{prop:conjugation-twist-cocy} are satisfied,
\begin{enumerate}
\item Group $G$ is abelian;
\item The codimension $l(g)$ is 1 or 2. $C(g)$ acts on $N^g$ by
isometry. When $l(g)=1,2$, isometry group of $N^g$ is abelian.
\end{enumerate}

At the end of this subsection, we give an example of the twisted
cyclic cocycle in a very special case.
\begin{ex}
Let $V$ be $\reals$, and let $G=\integers/2\integers=\{id, e\}$
act on $\reals$ by $e: x\mapsto -x$. In this case, $\Omega_e\in
C^1(A, A_e)$ is defined to be \[
\Omega_e(f)(x)=\frac{f(x)-f(-x)}{2x}.
\]
The cohomology $HH^\bullet (A, A_e)$ is computed to be
\[
HH^\bullet (A,  A_e)=\left\{\begin{array}{ll}0&\bullet \ne 1\\
\reals& \bullet =1\end{array}\right. ,
\]
where $HH^1(A, A_e)$ is generated by $\Omega_e$.
\end{ex}
\subsection{Inverse map of $L$}

We use the twisted cocycle constructed in the previous step to
obtain an inverse map of $L$
\[
T:  \Gamma^\infty(\bigoplus_{g\in G}\wedge^{\bullet
-l(g)}TV^g\otimes
\wedge^{l(g)}N^g)^G\longrightarrow C^\bullet
(Poly(V)\rtimes G, Poly(V)\rtimes G).
\]

We define $T$ as a composition of two maps $T_1$ and $T_2$. The map $T_1$ is defined as
\[
T_1: \Gamma^\infty(\bigoplus_{g\in G}\wedge^{\bullet
-l(g)}TV^g\otimes
\wedge^{l(g)}N^g)^{G}\longrightarrow C^\bullet (A,
Poly(V)\rtimes G)^G.
\]
The map $T_2$ is the standard map constructed in the proof of the
Eilenberg-Zilber theorem:
\[
T_2:C^\bullet (A,Poly(V)\rtimes G)^G\subset C^0(G,
C^\bullet (A, Poly(V)\rtimes G))\longrightarrow C^\bullet (A\rtimes
G, Poly(V)\rtimes G).
\]

\noindent{\bf Step I:} the map $T_1$ is a sum of the maps
\[
T_1^g:\Gamma^\infty(\wedge^{\bullet
-l(g)}TV^g\otimes
\wedge^{l(g)}N^g)\longrightarrow C^\bullet (A,
A_g),
\]
which are defined as follows.

Given $\xi\in \Gamma^\infty(\wedge^{\bullet
-l(g)}TV^g\otimes \wedge^{l(g)}N^g)$, we write
$\xi$ to be $X\otimes \Lambda_g$, where $X\in
\Gamma^\infty(\wedge^{\bullet -l(g)}TV^g)$ and
$\Lambda_g\in \Gamma^\infty(\wedge^{l(g)}N^g)$ is
defined be
\begin{equation}\label{eq:lambda}
\frac{\partial}{\partial x^{n-l(g)+1}}\wedge\cdots \wedge
\frac{\partial}{\partial x^{n}}.
\end{equation}

We define
\[
T_1^g(\xi)=X\sharp \Omega_g,\ \ \text{for any } \xi\in
\Gamma^\infty(\wedge^{k -l(g)}TV^g\otimes \wedge^{l(g)}N^g),
\]
where we view $X$ as a multidifferential operator on
$Poly(V)$ for $X\in \wedge^{k-l(g)}TV^g$ and $X\sharp \Omega_g(f_1,
\dots , f_k)$ equal to
\[
X(f_1, \dots , f_{k-l(g)})\Omega_g(f_{k-l(g )+1}, \dots , f_{k}).
\]


When $C(g)$ action on $N^g$ is diagonizable, as is explained in
Proposition \ref{prop:conjugation-twist-cocy}, we have
\[
T_1^g(h(\xi))=h(T_1^{h^{-1}gh}(\xi)),\ \ \ \forall h\in G.
\]
We define $T_1$ to be the sum of $T_1^g$. The restriction of
$T_1$ to the $G$-invariant sections gives the desired map
\[
T_1: \Gamma^\infty(\bigoplus_{g\in G}\wedge^{\bullet
-l(g)}TV^g\otimes
\wedge^{l(g)}N^g)^{G}\longrightarrow C^\bullet (A,
Poly(V)\rtimes G)^G.
\]
When $C(g)$ action on $N^g$ is not diagonalizable, then the image of
the above $T_1$ map are not always $G$-invariant. Therefore, we need
to replace $T_1$ by
\[
\frac{1}{|G|}\sum_h h(T_1(\xi)).
\]

\noindent{\bf Step II:} we now explain the construction of $T_2$, which is standard in
the Eilenberg-Zilber theorem. Given $\Phi\in C^k(A, Poly(V)\rtimes
G)$,
\begin{equation}\label{eq:t-2}
T_2(\Phi)(a_1U_{g_1}, \dots , a_kU_{g_k})=\Phi(a_1, \dots ,
g_1\cdots g_{k-1}(a_k))U_{g_1\cdots g_k}.
\end{equation}

By this Lemma \ref{lem:linear-twisted}, we have the following
proposition for the map $L_2$.
\begin{prop}
\label{prop:L2}Given $\xi\in
\Gamma^\infty(\oplus_g\wedge^{k-l(g)} TV^g\otimes
\wedge^{l(g)}N^g)^G$, we write $\xi=\sum_g
X_g\otimes \Lambda_g$, where $\Lambda_g$ is defined
 as in Equation (\ref{eq:lambda}).

The composition map $L_2 \circ T_1$ satisfies
\[
L_2(T_1(\xi))=\sum_g X_g\otimes \Lambda_g=\xi.
\]
\end{prop}
\begin{proof} As $L_2$ is $G$ equivariant and $\xi$ is $G$ invariant. We compute $L_2(T_1(\xi))(x)$ as follows: $L_2(T_1(\xi))(x)=$
\[
\begin{split}
=&\sum_{i_1, \dots , i_k} T_1(\xi)((x_1-x)^{i_1}
\cdots (x_k-x)^{i_k})\frac{\partial}{\partial
x^{i_1}}\wedge \cdots
\wedge \frac{\partial }{\partial x^{i_k}}\\
=&\sum_g\sum_{i_1, \dots , i_k} T_1(X_g\otimes
\Lambda_g)((x_1-x)^{i_1}\cdots
(x_k-x)^{i_k})\frac{\partial}{\partial x^{i_1}}\wedge \cdots
\wedge \frac{\partial }{\partial x^{i_k}}\\
=&\sum_g \sum_{i_1, \dots , i_k}X_g((x_1-x)^{i_1}, \dots
,
(x_{k-l(g)}-x)^{i_{k-l(g)}})\\
&\times\Omega_g((x_{k-l(g)+1}-x)^{i_{k-l(g)+1}},
\dots , (x_{k}-x)^{i_k})\frac{\partial}{\partial x^{i_1}}\wedge
\cdots \wedge \frac{\partial }{\partial x^{i_k}}\\
=&\sum_{g} \sum_{i_1, \dots , i_k}X_g ((x_1-x)^{i_1},
\dots , (x_{k-l(g)}-x)^{i_{k-l(g)}})
\frac{\partial}{\partial x^{i_1}}\wedge\cdots \wedge
\frac{\partial}{\partial
x^{i_{k-l(g)}}}\\
&\otimes\Omega_g((x_{k-l(g)+1}-x)^{i_{k-l(g)+1}},
\dots , (x_{k}-x)^{i_k})\frac{\partial}{\partial
x^{i_{k-l(g)+1}}}\wedge \cdots \wedge
\frac{\partial}{\partial x^{i_{k}}}\\
=&\sum_g X_g\otimes \Lambda_g.
\end{split}
\]
In the last equality, we used Lemma \ref{lem:linear-twisted} which shows $\Lambda_g$ and $\Omega_g$ have the same values on linear functions.
\end{proof}

We define $T=T_2\circ T_1$, and have the following theorem.

\begin{thm}
\label{thm:quasi-inverse} The map $T$ is a quasi-isomorphism. In
particular, $L\circ T=id$.
\end{thm}

\begin{proof}
We notice that $L_1(T_2)=id$ on $C^\bullet (A;A\rtimes
G)^G$, and therefore have
\[
L\circ T(\xi)=L_3(L_2(L_1(T_2(T_1(\xi)))))=L_3(L_2(T_1(\xi))),
\]
which is equal to $\xi$ by Proposition \ref{prop:L2}
\end{proof}
\subsection{The case of a smooth manifold}
In this subsection, we discuss the extension of the construction of
$T$ to general manifolds, which is again a composition of $T_2$ and $T_1$.

The map $T_2$ generalizes to the manifold case directly
because it is purely algebraic. The same formula as Equation
(\ref{eq:t-2}) defines a map
\[
T_2:\Big(C^\bullet\big(C^\infty(M), C^\infty(M)\rtimes
G\big)\Big)^G\longrightarrow C^\bullet\big(C^\infty(M)\rtimes G,
C^\infty(M)\rtimes G\big).
\]

To define $T_1$, we consider a $g$-invariant tubular
neighborhood $\calm^g$ of $M^g$. The neighborhood $\calm^g$ is a fiber bundle over
$M^g$, and we fix a $G$-invariant Ehresmann connection on $\calm^g$.
Furthermore, we choose a cut-off function $\rho_g$ on $\calm^g$
which is equal to 1 on a $g$-invariant neighborhood of $M^g$ and
vanishes outside $\calm^g$. Given a section $\xi_g =X_g\otimes
\Lambda_g$ of $\wedge^{\bullet -l(g)}TM^g \otimes \wedge^{l(g)}N^g$
to $\calm^g$, we use the Ehresmann connection to extend $X_g$ to a
multi-vector field $\tilde{X_g}$ on $\calm^g$, and define $\Omega_g$
a linear map on $C^\infty(\calm^g)$ by the same formula as Equation
(\ref{eq:twisted-cocycle}).

We define a cochain map
\[
T_1: \Gamma^\infty(\oplus _g \wedge^{k-l(g)}TM^g\otimes
\wedge^{l(g)}N^g)^G\longrightarrow C^k(C^\infty(M);
C^\infty(M)\rtimes G)^G
\]
by
\[
T_1(\xi)(f_1, \dots , f_k)=\sum_g
\rho_g\tilde{X}_g(f_1, \dots ,
f_{k-l(g)})\Omega_g(f_{k-l(g)+1}, \dots ,
f_{k})U_g,
\]
for $\xi=\sum_g X_g\otimes \Lambda_g$.

Again, we can easily compute that $L\circ T=id$. Therefore, since
$L$ is a quasi-isomorphism proved in \cite{nppt}, $T$ is also a
quasi-isomorphism.

\begin{rmk}
\label{rmk:canon} The definition of $T_1$ depends on a
choice of the normal bundle $N^g$, the $G$-invariant Ehresmann
connection, and the cut-off function. Therefore, $T$ is not a
canonical map. However, we notice that at any $x\in M^g$, inside
$T_xM$, there is a canonical complementary subspace to
$T_xM^g\subset T_xM$ determined by the representation of $g$ on
$T_xM$ independent of the choices of the metrics. Therefore, it is
easy to check ( {\it c.f.} Proposition \ref{prop:L2}) that when
restricted to the $\infty-$jets of $M^g$, the map $T$ is independent of all
the choices.
\end{rmk}

\section{Gerstenhaber bracket}

In this section, we compute the Gerstenhaber bracket on the
Hochschild cohomology of the algebra $C^\infty(M)\rtimes G$. Since
all the computations and constructions are local with respect to the
orbifold $M/G$, it is sufficient to work out everything locally on a
vector space. Because the Gerstenhaber bracket on
$C^\bullet(Poly(V)\rtimes G, Poly(V)\rtimes G)$ is the commutator of
the pre-Lie product, we mainly work on understanding the
pre-Lie product, and state the results for the Gerstenhaber bracket.

\subsection{Geometric description}
Let $\xi\in \Gamma^\infty\big(\oplus_{\alpha\in
G}\wedge^{k-l(\alpha)}TV^\alpha\otimes
\wedge^{l(\alpha)}N^\alpha\big)^G$, and $\eta \in
\Gamma^\infty\big(\oplus_{\beta\in G}\wedge^{l-l(\beta)}
TV^\beta\otimes \wedge^{l(\beta)}N^\beta\big)^G$. We compute the
pre-Lie product between $T(\xi)$ and $T(\eta)$ by $L(T(\xi)\circ
T(\eta))$. We write $\xi$ as the sum of
$\xi_\alpha=X_{\alpha}\otimes \Lambda_\alpha$, and $\eta$ as the sum
of $\eta_\beta=Y_\beta\otimes \Lambda_\beta$, with $X_\alpha\in
\Gamma^\infty(\wedge^{k-l(\alpha)}TV^\alpha),\ Y_\beta\in
\Gamma^\infty(\wedge^{l-l(\beta)}TV^\beta)$, $\Lambda_\alpha\in
\Gamma^\infty(\wedge^{l(\alpha)}N^\alpha)$, $\Lambda_\beta\in
\Gamma^\infty(\wedge^{l(\beta)}N^\beta)$. We compute the
Gerstenhaber bracket between $T(\xi)$ and $ T(\eta)$ by $ L([T(\xi),
T(\eta)])$ using the information of the quasi-isomorphisms $L$ and
$T$ introduced in Section 2 and 3.\\

\noindent{\bf Step I:} We compute $L_1(T(\xi)\circ T(\eta))\in
C^{k+l-1}(Poly(V), Poly(V)\rtimes G)$ first. As $\xi$
and $\eta$ are both $G$-invariant, the cocycles $T_1(\xi), T_1(\eta)\in
C^\bullet(Poly(V),\ Poly(V)\rtimes G)$ are also $G$-invariant, and
therefore the averaging in the definition of $L_1$ is not necessary. Our computations yield:
\[
\begin{split}
&L_1\big(T(X)\circ T(Y)\big)(f_1, \dots , f_{k+l-1})\\
=&\sum
\limits_s(-1)^{(s-1)(l-1)}T(\xi)\big(f_1, \dots , f_s, T(\eta)(f_{s+1}, \dots , f_{s+l}), f_{s+l+1}, \dots , f_{k+l-1}\big)\\
=&\sum\limits_s(-1)^{(s-1)(l-1)}T(\sum\limits_{\alpha}\xi_\alpha)\big(f_1, \dots , f_s,
T(\sum\limits_\beta \eta_\beta)(f_{s+1}, \dots , f_{s+l}), f_{s+l+1}, f_{k+l-1}\big)\\
=&\sum\limits_\alpha\sum\limits_\beta \sum\limits_s (-1)^{(s-1)(l-1)}T(\xi_\alpha)\big(f_1, \dots , f_s, T(\eta_\beta)(f_{s+1}, \dots , f_{s+l}),\\
&\qquad \qquad f_{s+l+1}, \dots , f_{k+l-1}\big)\\
=&\sum\limits_{g}\sum_{\alpha\beta=g}\sum\limits_s
(-1)^{(s-1)(l-1)}T_1^\alpha(\xi_\alpha)\big(f_1, \dots , f_s,
T_1^\beta(\eta_\beta)(f_{s+1}, \dots , f_{s+l}),\\
&\qquad \qquad \beta(f_{s+l+1}), \dots , \beta(f_{k+l-1})\big)U_{g}
\end{split}
\]

\noindent{\bf Step II:} We compute $L_2^g (L_1(T(\xi)\circ
T(\eta)))(x)$ by applying $L_2^g$ to the computation in the previous step: {\small
\begin{equation}
\label{eq:prelie}
\begin{split}
\frac{1}{|G|^2}&\sum \limits_{i_1, \dots ,
i_{k+l-1}}\sum\limits_{\alpha\beta=g, h_1, h_2\in G}\sum\limits_s
(-1)^{(s-1)(l-1)}h_1(T_1^\alpha(\xi_\alpha))\big((x_1-x)^{i_1},
\dots ,
(x_s-x)^{i_s},\\
&\qquad \qquad h_2(T_1^\beta(\eta_\beta))((x_{s+1}-x)^{i_{s+1}},
\dots, (x_{s+l}-x)^{i_{x+l}}), \beta((x_{s+l+1}-x)^{i_{s+l+1}}), \dots ,\\
&\qquad \qquad \beta((x_{k+l-1}-x)^{i_{k+l-1}})\big)
\frac{\partial}{\partial x^{i_1}}\wedge\cdots\wedge
\frac{\partial}{\partial x^{i_{k+l-1}}}.
\end{split}
\end{equation}
}

We look at the term $T_1^\beta(\eta_\beta)((x_{s+1}-x)^{i_{s+1}},
\dots , (x_{s+l}-x)^{i_{s+l}}))$. Using the expression
$\eta_\beta=X_\beta\otimes \Lambda_\beta$, we have
\[
\begin{split}
&T_1^\beta(\eta_\beta)\big((x_{s+1}-x)^{i_{s+1}}, \dots ,
(x_{s+l}-x)^{i_{s+l}})\big)\\
=&X_\beta\sharp \Omega_\beta \big((x_{s+1}-x)^{i_{s+1}}, \dots ,
(x_{s+l}-x)^{i_{s+l}})\big)\\
=&X_\beta\big((x_{s+1}-x)^{i_{s+1}},
\dots ,(x_{s+l-l(\beta)}-x)^{i_{s+l-l(\beta)}}\big)\\
\times&\Lambda_\beta\big((x_{s+l-l(\beta)}-x)^{i_{s+l-l(\beta)+1}},\dots
,
(x_{s+l}-x)^{i_{s+l}}\big)\\
=&\eta_\beta\big((x_{s+1}-x)^{i_{s+1}}, \dots,
(x_{s+l}-x)^{i_{s+l}}\big).
\end{split}
\]
In the second equality of the above equation, we used Lemma
\ref{lem:linear-twisted} which stated that $\Omega_\beta$ and $\Lambda_\beta$
agree on linear functions.

Substituting the above expression of the $T_1^\beta(\eta_\beta)$
into Equation (\ref{eq:prelie}), as $\eta$ is $G$ invariant, we have

\[
\begin{split}
&L^g_2\big(T(\xi_\alpha)\circ T(\eta_\beta)\big)\\
=&\frac{1}{|G|}\sum_h\sum \limits_{i_1, \dots ,
i_{k+l-1}}\sum\limits_s
(-1)^{(s-1)(l-1)}h(T_1^\alpha(\xi_\alpha))\Big((x_1-x)^{i_1}, \dots
,
(x_s-x)^{i_s},\\
&\qquad \qquad \eta_\beta\big((x_{s+1}-x)^{i_{s+1}},
\dots,(x_{s+l}-x)^{i_{s+l}}\big), \beta\big((x_{s+l+1}-x)^{i_{s+l+1}}\big), \dots ,\\
&\qquad \qquad
\beta\big((x_{k+l-1}-x)^{i_{k+l-1}}\big)\Big)\frac{\partial}{\partial
x^{i_1}}\wedge\cdots\wedge \frac{\partial}{\partial x^{i_{k+l-1}}}.
\end{split}
\]

We discuss several properties of $L_2^g(T(\xi_\alpha)\circ
T(\eta_\beta))$.
\begin{lem}
\label{lem:intersection}The restriction of $L(T(\xi_\alpha)\circ
T(\eta_\beta))$ to $V^\alpha\cap V^\beta$ is
\begin{equation}\label{eq:circ-twist}
\begin{split}
&\sum \limits_{i_1, \dots , i_{k+l-1}}\sum\limits_s
(-1)^{(s-1)(l-1)}\xi_\alpha\Big((x_1-x)^{i_1}, \dots ,
(x_s-x)^{i_s},\\
&\eta_\beta\big((x_{s+1}-x\big)^{i_{s+1}}, \dots,(x_{s+l}-x)^{i_{s+l}}\big), \beta\big((x_{s+l+1}-x)^{i_{s+l+1}}\big),
\dots ,\\
&\beta\big((x_{k+l-1}-x)^{i_{k+l-1}}\big)\Big)\frac{\partial}{\partial
x^{i_1}}\wedge\cdots\wedge \frac{\partial}{\partial x^{i_{k+l-1}}}.
\end{split}
\end{equation}
\end{lem}
\noindent{\em Proof:} Because $T_1(\xi_\alpha)\in C^k(Poly(V);
Poly(V)_\alpha)$ and $T_1(\eta_\beta)\in C^l(Poly(V);Poly(V)_\beta)$,
the cocycle $T_1(\xi_\alpha)\circ T_1(\eta_\beta)$ is in $C^{k+l-1}(Poly(V);
Poly(V)_{\alpha\beta})$. Therefore, the restriction of
$L(T(\xi_\alpha)\circ T(\eta_\beta))$ to $V^\alpha \cap V^\beta$ is
same to the restriction of $L_2^{\alpha\beta}(T_1(\xi_\alpha)\circ
T_1(\eta_\beta)))$ to $V^\alpha\cap V^\beta$.

We observe that restriction of $\Omega_\alpha(f_1, \dots ,
f_{l(\alpha)})$ (and $\Omega_\beta(\cdots)$) to $V^\alpha$ (and
$V^\beta$) is same to $\Lambda_\alpha(f_1, \dots , f_{l(\alpha)})$
(and $\Lambda_\beta(\cdots)$) as we set all the variables along the
normal direction of $V^\alpha$ equal 0. Using this property and $G$
invariance of $\xi$, we have that the restriction of
\[
\begin{split}
&T_1(\xi_\alpha)\Big((x_1-x)^{i_1}, \dots , (x_s-x)^{i_s},
T_1(\eta_\beta)\big((x_{s+1}-x)^{i_{s+1}}, \dots  (x_{s+l}-x)^{i_{s+l}}\big),\\
&\beta\big((x_{s+l+1}-x)^{i_{s+l+1}}\big), \dots ,
\beta\big((x_{k+l-1}-x)^{i_{k+l-1}}\big)\Big)
\end{split}
\]
to $V^\alpha\cap V^\beta$ is same to
\[
\begin{split}
&\xi_\alpha\Big((x_1-x)^{i_1}, \dots , (x_s-x)^{i_s},
\eta_\beta\big((x_{s+1}-x)^{i_{s+1}}, \dots,  (x_{s+l}-x)^{i_{s+l}}\big),\\
&\beta\big((x_{s+l+1}-x)^{i_{s+l+1}}\big), \dots ,
\beta\big((x_{k+l-1}-x)^{i_{k+l-1}}\big)\Big).\ \ \ \ \ \ \ \Box
\end{split}
\]

Inspired by the results of Lemma \ref{lem:intersection}, we
introduce the following definition.
\begin{dfn}
\label{dfn:tw-prelie}Let $V$ be a vector space with a linear
endomorphism $\gamma$. For all $\xi\in \wedge^k TV$, $\eta \in
\wedge^l TV$ , the $\gamma$-twisted pre-Lie product $\xi\circ_\gamma
\eta$ is defined to be
\[
\begin{split}
&\sum\limits_{i_1, \dots , i_{k+l-1}}\sum\limits_s
(-1)^{(s-1)(l-1)}\xi\Big((x_1-x)^{i_1}, \dots , (x_s-x)^{i_s},\\
&\eta\big((x_{s+1}-x)^{i_{s+1}}, \dots  (x_{s+l}-x)^{i_{s+l}})\big),
\gamma((x_{s+l+1}-x)^{i_{s+l+1}}), \dots ,\\
&\gamma((x_{k+l-1}-x)^{i_{k+l-1}}))\Big)(x)\frac{\partial}{\partial
x^{i_1}}\wedge \cdots \wedge \frac{\partial}{\partial
x^{i_{k+l-1}}}.
\end{split}
\]
\end{dfn}

\noindent{\bf Step III:} In the following, we look at the
$\alpha\beta$ component of the Gerstenhaber bracket in the case when
$V^\alpha\cap V^\beta=V^{\alpha\beta}$.


\begin{lem}
\label{lem:bracket-y}Assume $V^\alpha\cap V^\beta=V^{\alpha\beta}$.
If $\xi_{\alpha}=X_\alpha \otimes \Lambda_\alpha$ contains
directions along $N^\beta$, then the
contribution of $\xi_\alpha$ and $\eta_\beta$ to the $\alpha\beta$ component of
$L(T(\xi_{\alpha})\circ T(\eta_{\beta}))$ vanishes.
\end{lem}
\begin{proof}
Firstly, as $V^{\alpha\beta}=V^\alpha\cap V^\beta$, we can apply
Lemma \ref{lem:linear-twisted} to compute $L(T(\xi_{\alpha})\circ
T(\eta_{\beta}))$ using the twisted pre-Lie product between
$\xi_\alpha$ and $\eta_\beta$.

Secondly, we look at the term $\eta_\beta(\cdots)$ in
$X_\alpha\circ_\alpha \eta_\beta$ at $g=\alpha\beta$. It contains
derivations along all the directions of $N^\beta$.

Thirdly, by Lemma \ref{lem:linear-twisted}, if the component of
$L(T(\xi_\alpha)\circ T(\eta_\beta))$ at $\frac{\partial}{\partial
x^{i_1}}\wedge \cdots \wedge \frac{\partial}{\partial
x^{i_{k+l-1}}}$ does not vanish, we must have at least one of
$(x_1-x)^{i_1}, \dots , (x_s-x)^{i_s}$, $
\eta_\beta\big((x_{s+1}-x\big)^{i_{s+1}},
\dots,(x_{s+l}-x)^{i_{s+l}}\big),
\beta\big((x_{s+l+1}-x)^{i_{s+l+1}}\big), \dots ,
\beta\big((x_{k+l-1}-x)^{i_{k+l-1}}\big)$ which is supported along
the $N^\beta$ direction because $X_\alpha\otimes \Lambda_\alpha$
contains derivations along $N^\beta$. Furthermore, we notice that
$\eta_\beta\big((x_{s+1}-x\big)^{i_{s+1}},
\dots,(x_{s+l}-x)^{i_{s+l}}\big)$ lies along the $V^\beta$
direction. This implies that at least one of $(x_1-x)^{i_1}, \dots ,
(x_s-x)^{i_s}, \beta\big((x_{s+l+1}-x)^{i_{s+l+1}}\big), \dots ,
\beta\big((x_{k+l-1}-x)^{i_{k+l-1}}\big)$ must be along the
$N^\beta$ direction. Noticing that the action of $\beta$ on $V$
preserves the decomposition $V=V^\beta\oplus N^\beta$, we conclude
that at least one of $(x_1-x)^{i_1}, \dots , (x_s-x)^{i_s},
(x_{s+l+1}-x)^{i_{s+l+1}}, \dots , (x_{k+l-1}-x)^{i_{k+l-1}}$ is
along $N^\beta$.

Summarizing above observations, we see that $T(X_\alpha\otimes
\Lambda_\alpha)\circ T(Y_\beta\otimes\Lambda_\beta)$ contains too
many derivations along $N^\beta$ since in the above two parts
 we have a total of more than $l(\beta)$ many  derivations along
$N^\beta$ direction, which is of dimension $l(\beta)$. We conclude
that
\[
L(T(X_\alpha\otimes \Lambda_\alpha)\circ T(Y_\beta\otimes
\Lambda_\beta))= L(T(\xi_{\alpha})\circ T(\eta_{\beta}))=0.
\]
\end{proof}

Lemma \ref{lem:bracket-y} shows that to compute $\alpha\beta$
component of $L(T(\xi_\alpha)\circ T(\eta_\beta))$ we can assume
that $\xi_\alpha$ is contained in the $V^\beta$ direction. Therefore,
there is no need to consider the $\beta$ action on
$(x_{s+l+1}-x)^{i_{s+l+1}}, \cdots, (x_{k+l-1}-x)^{i_{k+l-1}}$ in  Equation
(\ref{eq:circ-twist}), which concerns $L(T(\xi_\alpha)\circ
T(\eta_\beta))$. In
this case, the twisted pre-Lie product between $\xi_\alpha$ and
$\eta_\beta$ is reduced to
\[
\begin{split}
&\sum\limits_{i_1, \dots , i_{k+l-1}}\sum\limits_s
(-1)^{(s-1)(l-1)}\xi_\alpha\Big((x_1-x)^{i_1}, \dots , (x_s-x)^{i_s},\\
&\qquad \qquad \eta_\beta\big((x_{s+1}-x)^{i_{s+1}}, \dots
(x_{s+l}-x)^{i_{s+l}}\big),\\
&\qquad \qquad \qquad \qquad (x_{s+l+1}-x)^{i_{s+l+1}}, \dots ,
(x_{k+l-1}-x)^{i_{k+l-1}}\Big)\frac{\partial}{\partial
x^{i_1}}\wedge \cdots \wedge \frac{\partial}{\partial
x^{i_{k+l-1}}}.
\end{split}
\]
This is the standard pre-Lie product $\xi_\alpha\circ \eta_\beta$
between $\xi_\alpha$ and $\eta_\beta$

We summarize the above computation into the following
theorem\footnote{Theorem similar to the one that follows was
stated in the first version of \cite{anno:preprint}, but the proof
there contained a crucial gap.}
\begin{thm}\label{thm:prelie}
Consider
\[
\begin{array}{c}
\xi=\sum_\alpha \xi_{\alpha}\in
\Gamma^\infty(\bigoplus_\alpha\wedge^{k-l(\alpha)}TV^\alpha
\otimes \wedge^{l(\alpha)}N^\alpha)^{G},\\
\eta=\sum_\beta \eta_{\beta}\in
\Gamma^\infty(\bigoplus_\beta\wedge^{l-l(\beta)}TV^\beta \otimes
\wedge^{l(\beta)}N^\beta)^{G}.
\end{array}
\]
If $V^\alpha\cap V^\beta=V^{\alpha\beta}$, then the contribution of $\xi_\alpha$ and $\eta_\beta$ to the $\alpha\beta$ component of
$L(T(\xi_{\alpha})\circ T(\eta_{\beta}))$ is
\[
pr^{\alpha\beta}(\xi_\alpha\circ \eta_\beta|_{V^{\alpha\beta}}),
\]
where $\circ$ in the above formula is the standard pre-Lie product
on $V$.
\end{thm}
\begin{proof}
Straight forward from Lemma \ref{lem:bracket-y}.
\end{proof}

We discuss a special but important case that $V^{\alpha}\cap
V^\beta=V^{\alpha\beta}$ holds true.  The following lemma is
well-known and we include its proof for readers' convenience.
\begin{lem}\label{lem:codim}
For $\alpha, \beta\in G$, the condition $V^\alpha+V^\beta=V$ is
equivalent to the equality that $l(\alpha)+l(\beta)=l(\alpha\beta)$.
And when $V^\alpha+V^\beta=V$,  $V^\alpha\cap
V^\beta=V^{\alpha\beta}$.
\end{lem}
\begin{proof}
As
$\dim(V^\alpha)+\dim(V^\beta)=\dim(V^\alpha+V^\beta)+\dim(V^\alpha\cap
V^\beta)$, we have that
\[
\begin{split}
&l(\alpha)+l(\beta)\\
=&\dim(V)-\dim(V^\alpha)+\dim(V)-\dim(V^\beta)\\
=&\dim(V)-\dim(V^\alpha+V^\beta)+\dim(V)-\dim(V^\alpha\cap V^\beta)\\
\geq&\dim(V)-\dim(V^\alpha+V^\beta)+\dim(V)-\dim(V^{\alpha\beta})\\
=&\dim(V)-\dim(V^\alpha+V^\beta)+l(\alpha\beta),
\end{split}
\]
where we used the fact that $V^\alpha\cap V^\beta\subset
V^{\alpha\beta}$.  Therefore, $l(\alpha)+l(\beta)=l(\alpha\beta)$
implies that $V=V^\alpha+V^\beta$.

On the other hand, assume that $V^\alpha+V^\beta=V$ and let $\langle\ ,\
\rangle$ be a $G$ invariant metric on $V$. For any $v\in
V^{\alpha\beta}$, we have $\alpha\beta(v)=v$, and accordingly
$\beta(v)=\alpha^{-1}(v)$, and $\beta(v)-v=\alpha^{-1}(v)-v$.
Furthermore, as the metric $\langle\ ,\ \rangle$  is $G$ invariant, we see that
$\beta(v)-v$ is orthogonal to $V^\beta$ with respect to the metric
$\langle\ ,\ \rangle$ and $\alpha^{-1}(v)-v$ is orthogonal to $V^\alpha$.
Therefore $\beta(v)-v=\alpha^{-1}(v)-v$ is orthogonal
to $V^\alpha+V^\beta$, which is equal to $V$ by the assumption. This
implies that $v$ must belong to $V^\alpha\cap V^\beta$, and we
have $V^\alpha\cap V^\beta=V^{\alpha\beta}$. This together with the
above equations implies that
\[
\begin{split}
l(\alpha)+l(\beta)=\dim(V)-\dim(V^\alpha+V^\beta)+\dim(V)-\dim(V^{\alpha\beta})=l(\alpha\beta).
\end{split}
\]
\end{proof}

\begin{ex}\label{ex:bracket-poisson}
We consider the set $S$ of all $g\in G$ such that $l(g)=2$. The set
$S$ is invariant under conjugation as $l(g)$ is invariant under
conjugate action. We consider $\pi=\sum_{g\in S}\pi_g\in
\Gamma^\infty(\bigoplus_{g\in S}\wedge^2N^g)^G$ and $\phi_{id}\in
\Gamma^\infty(\wedge^\bullet V)^G$, which is supported at the
identity component.

For $\pi$ and $\phi_{id}$, the conditions of Theorem \ref{thm:prelie}
are satisfied. Therefore, we have
\[
L(T(\pi)\circ T(\phi_{id}))=\sum_{g\in S}pr^g(\pi_g\circ \phi_{id}|_{V^g}),
\]
and
\[
L([T(\pi), T(\phi_{id})])=\sum_{g\in S}pr^g([\pi_g, \phi_{id}]|_{V^g}).
\]
\end{ex}

\subsection{Abelian group action}
In the following, we discuss the special case when $G$ is abelian.
Under this assumption, we obtain a more explicit description of the
twisted Schouten-Nijenhuis bracket. As $G$ is abelian, $G$ action on
$V$ is decomposed into a direct sum of 1 and 2 dimensional
irreducible representations of $G$, and we have global well defined
coordinate functions, which are eigenvectors for all $g$ action. In
particular, the conditions of Proposition
\ref{prop:conjugation-twist-cocy} are satisfied and the map $T_1$ in
Section 3.2 without averaging is already $G$ equivariant.

\begin{rmk}
Discussions in this subsection can be extended to the case where all elements in $(\alpha)$ (the set of elements in the same conjugacy class as $\alpha$) commute with all elements in $(\beta)$.
\end{rmk}

\begin{lem}
\label{lem:complete-intersection} Let
\[
\xi=\sum_{\alpha} X_\alpha\otimes \Lambda_\alpha\in
\Gamma^\infty(\bigoplus_\alpha \wedge ^{k-l(\alpha)}TV^\alpha\otimes
\wedge^{l(\alpha)}N^\alpha)^{G}
\]
and
\[
\eta=\sum_{\beta} Y_\beta\otimes \Lambda_\beta \in
\Gamma^\infty(\bigoplus_\beta \wedge^{l-l(\beta)} TV^\beta\otimes
\wedge^{l(\beta)}N^\beta)^G,
\]
if $V^\alpha+V^\beta\ne V$ for all $\alpha, \beta$, then
\[
L(T(\xi)\circ T(\eta))=0.
\]
\end{lem}
\begin{proof} Following the computations similar to those in  Equation
(\ref{eq:prelie}), we have that $L(T(\xi_{(\alpha)})\circ
T(\eta_{(\beta)}))$ is equal to
\[
\begin{split}
&\sum \limits_{i_1, \dots ,
i_{k+l-1}}\sum\limits_{\alpha\beta=g}\sum\limits_s
(-1)^{(s-1)(l-1)}T_1^\alpha(\xi_\alpha)\Big((x_1-x)^{i_1}, \dots ,
(x_s-x)^{i_s},\\
&\qquad \qquad T_1^\beta(\eta_\beta)\big((x_{s+1}-x)^{i_{s+1}},
\dots,
(x_{s+l}-x)^{i_{x+l}}\big),\\
&\beta\big((x_{s+l+1}-x)^{i_{s+l+1}}\big), \dots ,
\beta\big((x_{k+l-1}-x)^{i_{k+l-1}}\big)\Big)\frac{\partial}{\partial
x^{i_1}}\wedge\cdots\wedge \frac{\partial}{\partial x^{i_{k+l-1}}}.
\end{split}
\]

Since $V^\alpha+V^\beta\ne V$, then its normal directions
$N^\perp=N^\alpha\cap N^\beta$ are nontrivial. We observe that $T(\xi_\alpha)$
contains all the derivations\footnote{Rigorously speaking,
$\Omega_\alpha, \Omega_\beta$ are not derivations. Here, we use the
word ``derivation" loosely, because they behave like derivations on
linear functions.} along $N^\alpha$, and $T(\eta_\beta)$ contains
all the derivations along $N^\beta$. As $\eta_\beta$ is a section of
$\wedge^{l-l(\beta)}TV^\beta\otimes \wedge^{l(\beta)}N^\beta$,
$T_1^\beta(\eta_\beta)((x_{s+1}-x)^{i_{s+1}}, \dots ,
(x_{s+l}-x)^{i_{s+l}}))$ is independent of $N^\perp$. Therefore, we
see that to have a nonzero contribution in $L(T(\xi_\alpha)\circ
T(\eta_\beta))$, the coordinates $x^{i_1}, \cdots, x^{i_{s+l}},
\beta(x^{i_{s+l+1}}), \cdots$, $ \beta(x^{i_{k+l-1}})$ must contain two
copies of the variables along $N^\perp$ and one copy of the
variables along $N^\beta/N^\perp$. However, because $i_1, \dots ,
i_{k+l-1}$ are distinguished and $N^\beta$ is $\beta$-invariant, the coordinates
$x^{i_1}, \dots , x^{s+k-1},$ $ \beta(x^{s+k}),\dots ,
\beta(x^{i_{k+l-1}})$ have at most one copy of the variables along
the $N^\beta$ direction. There are not enough variables along the
$N^\perp$ direction. This implies that $L(T(\xi_{(\alpha)})\circ
T(\eta_{(\beta)}))$ vanishes.
\end{proof}

By Lemma \ref{lem:complete-intersection}, we are reduced to considering the
Gerstenhaber bracket in the case when $V^\alpha+V^\beta=V$, which by
Lemma \ref{lem:codim} implies $V^\alpha\cap V^\beta=V^{\alpha\beta}$.
Therefore, we can use Theorem \ref{thm:prelie} to compute the Gerstenhaber bracket.
\begin{lem}
\label{lem:bracket-x}Let $G$ be an abelian group. If $l-l(\beta)\geq
2$, for $Y_\beta\in \Gamma^\infty(\wedge
^{l-l(\beta)}TN^\alpha)^{C(\beta)}$,
\[
L(T(\xi_{\alpha})\circ T(\eta_{\beta}))=0.
\]
\end{lem}
\begin{proof} Because of Lemma \ref{lem:complete-intersection} and Theorem \ref{thm:prelie}, we can drop the $\beta$ twist in Equation (\ref{eq:circ-twist}). At the component $g=\alpha\beta$,
we look at the number of derivations along the direction of
$N^\alpha$. The component $T(\xi_\alpha)$ contributes $l(\alpha)$ and
$T(Y_\beta\otimes \Lambda_\beta)$ contributes $l-l(\beta)$ number of derivations.
Therefore, $T(\xi_\alpha)\circ T(Y_\beta\otimes \Lambda_\beta)$
contains at least the following number of derivations along
$N^\alpha$:
\[
l(\alpha)+(l-l(\beta))-1\geq l(\alpha)+2-1\geq l(\alpha)+1.
\]
This implies the statement of this lemma, because
$dim(N^\alpha)=l(\alpha)$.
\end{proof}

In the case when $G$ is abelian, the expressions of the pre-Lie
product and Gerstenhaber bracket in Theorem \ref{thm:prelie} can be
simplified.
\begin{thm}
\label{thm:main}Let
\[
\xi=\sum_{\alpha\in (\alpha)}X_\alpha \otimes \Lambda_\alpha\in
\Gamma^\infty(\wedge^{k-l(\alpha)}TV^\alpha\otimes
\wedge^{l(\alpha)}N^\alpha)^{C(\alpha)},
\]
and
\[
\eta=\sum_{\beta\in (\beta)}Y_\beta\otimes \Lambda_\beta \in
\Gamma^\infty(\wedge^{l-l(\beta)}TV^\beta\otimes
\wedge^{l(\beta)}N^\beta)^{C(\beta)}.
\]
Then the component of $ L(T(\xi)\circ T(\eta))$ in
\[
\Gamma^\infty(\wedge^{k+l-l(\alpha)-1}TV^{\alpha\beta}\otimes
\wedge^{l(\alpha)}N^{\alpha\beta})^{C(\alpha\beta)}
\]
is computed as follows.
\begin{enumerate}
\item  When $V^\alpha+V^\beta\ne V$ for all $\alpha, \beta$, then $L(T(\xi)\circ T(\eta))=0$.
\item  When $V^\alpha+V^\beta=V$, as $G$ is abelian, we write $V=V^{\alpha\beta}\oplus
N^\alpha\oplus N^\beta$, where $V^{\alpha\beta}$ is the invariant
subspace of $\alpha\beta$, $N^\alpha$ is the subspace orthogonal
to $V^\alpha$, and $N^\beta$ is the subspace orthogonal to
$V^\beta$. In this case, $V^\alpha=V^{\alpha\beta}\oplus N^\beta$,
and $V^\beta=V^{\alpha\beta} \oplus N^\alpha$. According to this
decomposition, we write
\[
\begin{split}
&X_\alpha=\sum X_\alpha^{i_1\cdots i_s, p_1\cdots p_{k-l(\alpha)-s}}\frac{\partial}{\partial x^{i_1}}\wedge \cdots \frac{\partial}{\partial x^{i_s}}\wedge \frac{\partial}{\partial x^{p_1}}\wedge \cdots \wedge \frac{\partial}{\partial x^{p_{k-l(\alpha)-s}}}\\
&\in \Gamma^\infty\big(\wedge^{s}TV^{\alpha\beta}\big)\otimes \Gamma^\infty\big( \wedge^{k-l(\alpha)-s}TN^\beta\big)\\
&Y_\beta=\sum Y_\beta^{j_1\cdots j_t, q_1\cdots q_{l-l(\beta)-t}}\frac{\partial}{\partial x^{j_1}}\wedge \cdots \frac{\partial}{\partial x^{j_t}}\wedge \frac{\partial}{\partial x^{q_1}}\wedge \cdots \wedge \frac{\partial}{\partial x^{q_{l-l(\beta)-t}}}\\
&\in \Gamma^\infty\big(\wedge^{t}TV^{\alpha\beta}\big)\otimes
\Gamma^\infty \big(\wedge^{k-l(\beta)-t}TN^\alpha \big).
\end{split}
\]

The component of $L(T(\xi_{(\alpha)})\circ T(\eta_{(\beta)}))$ in
\[
\Gamma^\infty(\wedge^{k+l-l(g)-1}TV^{g}\otimes
\wedge^{l(g)}(N^\alpha\oplus N^\beta))^{C(g)}
\]
is
computed to be
\[
\begin{split}
&\sum\limits_{\tiny \begin{array}{c}g=\lambda\mu\\ \lambda\in (\alpha), \mu\in
(\beta)\\ l(g)=l(\lambda)+l(\mu)\end{array}} \sum\limits_{i_1, \dots , i_{k-l(\lambda)}, j_1, \dots ,
j_{l-l(\mu)}}
(-1)^{(z-1)(l-1)+(k-z)l+(l-l(\mu))l(\lambda)}X_\lambda^{i_1\cdots\widehat{i_z}
\cdots i_{k-l(\lambda)}}\\
&\frac{\partial}{\partial x^{i_z}}Y_\mu
^{j_1\cdots j_{l-l(\mu)}}
\frac{\partial}{\partial x^{i_1}}\wedge \cdots
\widehat{\frac{\partial}{\partial x^{i_z}}} \cdots \wedge
\frac{\partial}{\partial x^{i_{k-l(\lambda)}}}\wedge \frac{\partial
}
{\partial x^{j_1}}\wedge\cdots \wedge\frac{\partial}{\partial x^{j_{l-l(\mu)}}}\otimes \Lambda_\lambda\wedge\Lambda_\mu\\
&+\sum\limits_{i_1, \dots , i_{k-l(\lambda)}, j_1, \dots ,
j_{l-l(\mu)-1}, q_z}
(-1)^{(k-l(\lambda))(l-1)+(k-l(\mu))(l(\lambda)-1)}X_\lambda^{i_1\cdots
i_{k-l(\lambda)}}\\
&\frac{\partial }{\partial x^{q_z}}Y_\mu^{j_1,
\dots , j_{l-l(\mu)-1},
q_{z}}
\frac{\partial}{\partial x^{i_1}}\wedge \cdots \wedge
\frac{\partial}{\partial x^{i_{k-l(\lambda)}}}\wedge
\frac{\partial}{\partial x^{j_1}}\wedge\cdots
\frac{\partial}{\partial x^{j_{l-l(\mu)-1}}}\otimes
\Lambda_\lambda\wedge\Lambda_\mu.
\end{split}
\]
\end{enumerate}
\end{thm}
\begin{proof}
The first statement is a corollary of Lemma
\ref{lem:complete-intersection}. It remains to show the second
statement.

As $G$ is abelian, we can simultaneously diagonalize the $\alpha,
\beta$ action on $V$. Using the fact $V^\alpha+V^\beta=V$, we can write
$V=V^{\alpha\beta}\oplus N^\alpha\oplus N^\beta$ such that
$V^\alpha=V^{\alpha\beta}\oplus N^\beta$ and
$V^\beta=V^{\alpha\beta}\oplus N^\alpha$.

By Lemma \ref{lem:bracket-y}, we conclude that $X_\alpha$ must be
from $\Gamma^\infty(\wedge^{k-l(\alpha)}TV^{\alpha\beta})$ to have
nontrivial contribution in $L([T(\xi_{(\alpha)}),
T(\eta_{(\beta)})])$.

Similarly, by Lemma \ref{lem:bracket-x}, we conclude that to have nontrivial
contribution in $L([T(\xi_{(\alpha)}), T(\eta_{(\beta)})])$, the direction
$Y_\beta$ has to be from one of the following spaces:
\[
\begin{split}
(i)\ Y_\beta\in
\Gamma^\infty(\wedge^{l-l(\beta)}TV^{\alpha\beta})\ \ \ \ \ \
&(ii)\ Y_\beta \in
\Gamma^\infty(\wedge^{l-l(\beta)-1}TV^{\alpha\beta})\otimes
\Gamma^\infty(TN^\alpha).
\end{split}
\]

We apply Theorem \ref{thm:prelie} to compute the pre-Lie
product. Because both $V^{\alpha\beta}$ and $N^\alpha$ are
subspaces of $V^\beta$ which is the fixed point set of $\beta$, we
can drop the $\beta$ twist of the pre-Lie product. Therefore, we
are left with the standard pre-Lie product.

When $Y_\beta\in
\Gamma^\infty(\wedge^{l-l(\beta)}TV^{\alpha\beta})$, we have
$L([T(\xi_\alpha), T(\eta_\beta)])=$
\[
\begin{split}
&\sum\limits_{i_1, \dots , i_{k-l(\alpha)}, j_1, \dots ,
j_{l-l(\beta)}}(-1)^{(z-1)(l-1)+(k-z)l+(l-l(\beta))l(\alpha)}
X_\alpha^{i_1\cdots\widehat{i_z}\cdots i_{k-l(\alpha)}}\\
&\frac{\partial}{\partial x^{i_z}}Y_\beta ^{j_1\cdots j_{l-l(\beta)}}
\frac{\partial}{\partial x^{i_1}}\wedge \cdots
\widehat{\frac{\partial}{\partial x^{i_z}}}\cdots \wedge
\frac{\partial}{\partial x^{i_{k-l(\alpha)}}}\wedge \frac{\partial
}{\partial x^{j_1}}\wedge\cdots \wedge\frac{\partial}{\partial
x^{j_{l-l(\beta)}}}.
\end{split}
\]

When $Y_\beta\in
\Gamma^\infty(\wedge^{l-l(\beta)-1}TV^{\alpha\beta})\otimes
\Gamma^\infty(TN^\alpha)$, we compute $L([T(\xi_\alpha),
T(\eta_\beta)])=$
\[
\begin{split}
&\sum\limits_{i_1, \dots , i_{k-l(\alpha)}, j_1, \dots
,j_{l-l(\beta)-1},
q_z}(-1)^{(k-l(\alpha))(l-1)+(k-l(\beta))(l(\alpha)-1)}
X_\alpha^{i_1\cdots i_{k-l(\alpha)}}\\
&\frac{\partial }{\partial
q_z}Y_\beta^{j_1, \dots , j_{l-l(\beta)-1}, q_{z}}\frac{\partial}{\partial x^{i_1}}\wedge \cdots \wedge
\frac{\partial}{\partial x^{i_{k-l(\alpha)}}}\wedge
\frac{\partial}{\partial x^{j_1}}\wedge\cdots
\frac{\partial}{\partial x^{j_{l-l(\beta)-1}}}.
\end{split}
\]
\end{proof}
\begin{cor}\label{cor:prelie}
Let $G$ be a abelian group acting a manifold $M$. The $(g)$
component of the twisted Schouten-Nijenhuis bracket is
\[
L([T(\xi), T(\eta)]) =pr^g(\sum_{\tiny
\begin{array}{c}g=\lambda\mu\\ \lambda\in (\alpha), \mu\in (\beta)\\
l(g)=l(\lambda)+l(\mu)\end{array}}[\tilde{\xi}_{\lambda},
\tilde{\eta}_\mu]|_{M^{\lambda\mu}}).
\]
\end{cor}

\section{Noncommutative Poisson structure and symplectic reflection algebras}
In this section, we want to find geometric expressions of all
possible noncommutative Poisson structures on $C^\infty(M)\rtimes
G$.

\subsection{Noncommutative Poisson structure}

For a noncommutative Poisson structure on $C^\infty(M)\rtimes G$, we
need to consider degree 2 Hochschild cohomology classes. Since there
is no contribution from $g$-fixed point submanifolds with $l(g)=1$
(Equation (\ref{eq:coh-even})), we have the following description of
degree 2 Hochschild cohomology of $C^\infty(M)\rtimes G$,
\[
HH^2(C^\infty(M)\rtimes G, C^\infty(M)\rtimes
G)=\Gamma^{\infty}(\wedge^2 TM)^G\oplus
\Gamma^\infty\big(\bigoplus_{g\in G, l(g)=2}\wedge^2N^g\big)^G,
\]
where the $G$-action on the second component is by conjugation. Inspired
by the above equation, we define the set $S$ of elements $g\in G$
such that the fixed point subspace of $g$ is of codimension 2. It is
easy to see that $S$ is closed under the conjugate action of $G$.

With the above discussion, we are ready to state the following
geometric description of noncommutative Poisson structures on
$C^\infty(M)\rtimes G$.

We assume that the manifold $M$ has a $G$-invariant symplectic
structure. As $G$ is finite, there always exists a $G$-invariant
compatible almost complex structure $J$ on $M$. Therefore, for any
$g$ in $G$, the fixed point manifold $M^g$ is a symplectic
submanifold with a compatible almost complex structure
\cite{gu-st:fixed-pt} and accordingly is of even dimension.
Furthermore, we have that the restriction of the $l(g)$-th wedge
power of the corresponding Poisson structure defines a global
section on $\wedge^{l(g)}N^g$, which is $G$ invariant.
\begin{thm}
\label{thm:poisson} Assume that $M$ is a symplectic
manifold\footnote{From the proof, we see that all we really need is
a $G$-invariant almost complex structure on $M$.}, with a $G$
symplectic action. Let $\pi$ be an element in
$\Gamma^\infty(\wedge^2TM)^{G}$, and $\sum_{g\in S} \Lambda_g$ be an
element in $(\oplus_{g\in S }\wedge^2 N^g)^G$. Then
$\pi+\sum_{g}\Lambda_g$ is a noncommutative Poisson structure on
$C^\infty(M)\rtimes G$ if and only if
\begin{enumerate}
\item On $M$, $[\pi, \pi]=0$;
\item For any $g\in S$, $pr^g([\pi, \widetilde{\Lambda}_g]|_{M^g})=0$.
\end{enumerate}

In the above, we  have assumed that there is no group element of $G$
except the identity that is a stabilizer of an open subset of $M$.
And in this case, we call the $G$ action is reduced.
\end{thm}
\begin{proof}We compute $L([T(\pi+\sum_{g\in S}\Lambda_g),
T(\pi+\sum_{g\in S}\Lambda_g)])$. It decomposes into the
sum of four terms
\[
\begin{split}
i)\ L([T(\pi), T(\pi)]),\ \ \ \ \ \ \ \ & ii)\
\sum\limits_{\alpha\in S}L([T(\pi), T(\Lambda_\alpha)]),\\ iii)\
\sum\limits_{\alpha\in S }L([T(\Lambda_\alpha), T(\pi)]),\ \ \ \ \
\ & iv)\ \sum\limits_{\alpha, \beta\in S}L([T(\Lambda_\alpha),
T(\Lambda_\beta)]).
\end{split}
\]

We compute the above i)-iv) terms separately.
\begin{enumerate}
\item $L([T(\pi), T(\pi)])$. On the identity component the
Gerstenhaber bracket corresponds to the standard
Schouten-Nijenhuis bracket. Therefore, we have
\[
L([T(\pi), T(\pi)])=[\pi, \pi],
\]
which is again on the identity component.
\item $\sum_{\alpha\in S}L([T(\pi), T(\Lambda_\alpha)])$. The identity component contributes $\pi$ and the $\alpha$ component contributes $\Lambda_\alpha$. We
use Example \ref{ex:bracket-poisson} to compute these terms,
\[
\sum_{\alpha\in S}L([T(\pi), T(\Lambda_\alpha)])=\sum_{\alpha\in
S}pr^\alpha([\pi, \widetilde{\Lambda}_\alpha]|_{M^\alpha}),
\]
where $[\pi, \widetilde{\Lambda}_\alpha]|_{M^\alpha}$ is on the
$\alpha$ component.
\item $\sum_{\alpha\in
S}L([T(\Lambda_\alpha), T(\pi)]))$. This case is similar to the previous
case. We apply Example \ref{ex:bracket-poisson} to compute the
terms
\[
\sum_{\alpha\in S}L([T(\Lambda_\alpha), T(\pi)])=\sum_{\alpha\in
S}pr^\alpha([\widetilde{\Lambda}_\alpha, \pi]_{M^\alpha}),
\]
where $[\widetilde{\Lambda}_\alpha, \pi]_{M^\alpha}$ is on the
$\alpha$ component.
\item $\sum_{\alpha, \beta\in
S}L([T(\Lambda_\alpha), T(\Lambda_\beta)])$. We need to use the
assumption that $M$ is a symplectic manifold. Accordingly,
$M^\alpha$ is a symplectic submanifold of codimension 2.
Furthermore, if we fix a $G$-invariant compatible almost complex
structure on $M$, $M^\alpha$ is an almost complex submanifold.
Hence, for each $x\in M^\alpha\cap M^\beta$ with $\alpha, \beta\in
S$, $T_xM^\alpha$ and $T_xM^\beta$ are almost complex subspaces of
$T_xM$. The sum $T_xM^\alpha+ T_xM^\beta$ is again an almost
complex subspace of $T_xM$ and satisfies one of the following properties.
\begin{enumerate}
\item $T_xM^\alpha+ T_xM^\beta=T_xM$,
\item $T_xM^\alpha=T_xM^\beta=T_xM^{\alpha\beta}$,
\item $T_xM^\alpha=T_xM^\beta$ with $T_xM^{\alpha\beta}=T_xM$.
\end{enumerate}
We discuss these three cases separately.

\noindent{(a)} If $T_xM^\alpha+ T_xM^\beta=T_xM$, then by Lemma
\ref{lem:codim}, $T_xM^{\alpha\beta}=T_xM^\alpha\cap
T_xM^\beta$. Therefore, $M^{\alpha\beta}$ is a codimension 4
submanifold near $x$. We notice that $[T(\Lambda_\alpha),
T(\Lambda_\beta)]$ is a 3-cocycle. By Equation (\ref{eq:coh-even})
 with $\bullet=3$, we see that any contribution to the degree 3 Hochschild cohomology
of $C^\infty(M)\rtimes G$ from a $g$-fixed point submanifold with
$l(g)=4$ is trivial. Hence, the $\alpha\beta$ component of
$\sum_{\alpha, \beta\in S}L([T(\Lambda_\alpha),
T(\Lambda_\beta)])(x)$ vanishes.

\noindent{(b)} If $T_xM^\alpha=T_xM^\beta=T_xM^{\alpha\beta}$, then
$M^\alpha=M^\beta=M^{\alpha\beta}=M^\alpha\cap M^\beta$ near $x$.
Therefore, we apply Lemma \ref{lem:intersection} to compute
$L(T(\Lambda_\alpha)\circ T(\Lambda_\beta))$ and
$L(T(\Lambda_\beta)\circ T(\Lambda_\alpha))$ by the twisted pre-Lie
product (Definition \ref{dfn:tw-prelie}). We observe that the normal
bundle $N^\alpha$ and $N^\beta$ are both symplectic orthogonal to
$M^\alpha=M^\beta$. Therefore, $N^\alpha$ is same as $N^\beta$. In
Definition \ref{dfn:tw-prelie}, since the $\beta$ action preserves
$N^\beta=N^\alpha$, to have a nontrivial outcome we need all 3 terms
$(x_s-x)^{i_s}$ for $s=1, \cdots, 3$ to be along the
$N^\alpha=N^\beta$ direction.
This is because $\Lambda_\alpha$ and $\Lambda_\beta$ both contains 2
derivations along $N^\alpha=N^\beta$. On the other hand,
$N^\alpha=N^\beta$ is only of 2 dimension. Any wedge product of 3
vectors along $N^\alpha$ vanishes. Therefore, the $\alpha\beta$
component of $\sum_{\alpha, \beta\in S}L([T(\Lambda_\alpha),
T(\Lambda_\beta)])(x)$ vanishes.

\noindent{(c)} If $T_xM^\alpha=T_xM^\beta$ but
$T_xM^{\alpha\beta}=T_xM$. This shows that $\alpha\beta$ acts on
$T_xM$ trivially and therefore there is a neighborhood of $x$ which
is totally fixed by $\alpha\beta$. Hence, by the assumption that $G$
action on $M$ is reduced, we know $\alpha\beta=1$. (As
$N^\alpha=N^\beta$ is of 2 dimension, the centralizer subgroups
$C(\alpha)$ and $C(\beta)$ action on $N^\alpha$, which are subgroups
of $U(1)$, can be diagonalized simultaneously. We choose the
eigenfunctions of this action to be the coordinate functions. See
Lemma \ref{prop:conjugation-twist-cocy}.) We notice that
$T(\Lambda_\alpha)$ and $T(\Lambda_\beta)$ both contains 2
derivations along $N^\alpha=N^\beta$ direction, it is not difficult
to check that $L([T(\Lambda_\alpha), T(\Lambda_\beta)])$ contains
wedge product of 3 vector fields along $N^\alpha=N^\beta$ direction.
Therefore, the identity component of $L([T(\Lambda_\alpha),
T(\Lambda_\beta)])(x)$ vanishes for the same reason as part (b).
\end{enumerate}
In summary, we have that
\[
\begin{split}
&L([T(\pi+\sum_{g\in S}\Lambda_g),
T(\pi+\sum\limits_{g}\Lambda_g)])\\
=&[\pi,\pi]+\sum\limits_{\alpha\in S}pr^\alpha([\pi,
\tilde{\Lambda}_\alpha]|_{M^\alpha})+\sum\limits_{\alpha\in
S}pr^\alpha([\tilde{\Lambda}_\alpha, \pi]|_{M^\alpha})\\
=&[\pi, \pi]+2\sum\limits_{\alpha\in S}pr^\alpha([\pi,
\tilde{\Lambda}_\alpha]|_{M^\alpha}).
\end{split}
\]
Therefore, if $\pi+\sum_{g\in S}\Lambda_g$ is a noncommutative
Poisson structure, then $[\pi,\pi]=0$ and $pr^\alpha([\pi,
\tilde{\Lambda}_\alpha]|_{M^\alpha})=0$ for any $\alpha\in S$.
\end{proof}

\begin{cor}
\label{cor:linear} Let $V$ be a real symplectic vector space with a
$G$ invariant linear symplectic form $\omega$. Assume that the $G$
action on $V$ is reduced as in Theorem \ref{thm:poisson}. Let $\pi$
be the corresponding Poisson structure of $\omega$. Then
$\kappa=\pi+\sum_{\alpha\in S}\Lambda_\alpha$ is a noncommutative
Poisson structure on $Poly(V)\rtimes G$ if and only if
$\Lambda_\alpha$ is constant on $V^\alpha$.
\end{cor}
\begin{proof}
Under the assumption of the corollary, the restriction of the
symplectic form $\omega$ to $N^\alpha$ for $\alpha\in S$ is a
symplectic two form. We denote the corresponding Poisson structure
on $N^\alpha$ by $\pi_\alpha$. Accordingly, we can write
$\Lambda_\alpha=f_\alpha\pi_\alpha$, where $f_\alpha$ is a
polynomial function on $V^\alpha$.

By Theorem \ref{thm:poisson},  $\pi+\sum_{\alpha\in
S}f_\alpha\pi_\alpha$ is a noncommutative Poisson structure if and
only if
\begin{enumerate}
\item $[\pi, \pi]=0$,
\item $[\pi, f_\alpha]|_{V^\alpha}=0$.
\end{enumerate}

Equation (1) is automatically satisfied because $\pi$ is Poisson.
Because $[\pi, \pi_\alpha]=0$, Equation (2) is reduced to
\[
pr^\alpha([\pi, f_\alpha\pi_\alpha]|_{V^\alpha})=pr^\alpha([\pi, f_\alpha]\wedge
\pi_\alpha|_{V^\alpha})=[\pi,f_\alpha]|_{V^\alpha}\wedge
\pi_\alpha=0.
\]
In the second equality, we have used the fact that as $f_\alpha$ is a function on $V^\alpha$, the bracket $[\pi, f_\alpha]$ is a vector field along $V^\alpha$.
Therefore, $\pi+\sum_{\alpha\in S}f_\alpha\pi_\alpha$ is a
noncommutative Poisson structure if and only if $[\pi,
f_\alpha]|_{V^\alpha}=[\pi, f_\alpha]=0$, for all $\alpha\in S$.  Because $\omega$ is symplectic, $[\pi,
f_\alpha]=0$ forces $f_\alpha$ to be a constant. Therefore,
$\Lambda_\alpha$ is also a constant on $V^\alpha$ for all $\alpha\in
S$.
\end{proof}
\subsection{Remarks on deformation quantizations}
It is known that the set of noncommutative Poisson structures on an algebra
$A$ is in one to one correspondence to the set of infinitesimal deformations
of $A$. It is natural to ask whether one can integrate the
infinitesimal deformation associated to a noncommutative Poisson structure to a real one. This question relates to
the idea of deformation quantization in mathematical physics. In
\cite{t1:def-gpd}, the second author introduced a notion of
deformation quantization of a noncommutative Poisson structure,
which we now recall.
\begin{dfn}
\label{dfn:quantization} A deformation quantization of a
noncommutative Poisson structure $\Pi$ on an associative algebra
$A$ is an associative product $\star_\hbar$ on $A[[\hbar]]$, such
that $f\star_\hbar g=\sum\limits_i \hbar^i C_i(f,g)$ for all $f,g\in
A$ satisfying the following conditions:
\begin{enumerate}
\item $C_0(f,g)=fg$;
\item The Hochschild cohomology class $[C_1]$ is equal to $\Pi$.
\end{enumerate}
\end{dfn}

It is natural to ask whether all the noncommutative Poisson
structures defined in Theorem \ref{thm:poisson} can be deformation
quantized. One special case has already been studied extensively, namely,
when the noncommutative Poisson structure comes from an $G$-invariant Poisson structure on $M$. (i.e. there are no components supported on fixed point submanifolds of codimension 2.)
The deformation quantizations of these types of noncommutative
Poisson structures on $C^\infty(M)\rtimes G$ were studied in
\cite{t1:def-gpd}, \cite{t2:thesis}, \cite{do-et:cohomology},
\cite{nppt}, etc. Another well-studied and well-known special case is the following
proposition, which is essentially due to Etingof and Ginzburg
\cite[Theorem 1.3]{eg:cherednik}.
\begin{prop}
\label{prop:cherednik}The noncommutative Poisson structure on a
symplectic vector space obtained in Corollary \ref{cor:linear} can
be deformation quantized.
\end{prop}
\begin{proof}In the following, we work with the field $\complex$,
because we will use the construction of the symplectic reflection
algebras in \cite{eg:cherednik} and Theorem 1.3 therein.
Everything extends to the field $\reals$, because Theorem 1.3 in
\cite{eg:cherednik} still holds in the real case. (The real group
algebra of a finite group is semisimple.)

In \cite{eg:cherednik}, a symplectic reflection algebra $H_{t,c}$
is introduced as
\[
TV\rtimes G/<x\otimes y-y\otimes x-\kappa(x,y)\in T^2V\oplus
\complex G>_{x,y\in V},
\]
where $(V, \omega)$ is a finite dimensional complex symplectic
vector space over $\complex$, and $TV$ is its tensor algebra, and
$T^2V$ is $V\otimes V$, and $\kappa$ is defined to be
\[
\kappa(x, y)=t\pi(x,y)+\sum_{\alpha\in S}c_\alpha\pi_\alpha
U_\alpha,
\]
a $G$-invariant section of $\wedge^2 V+\oplus_{\alpha\in S}
\wedge^2 N^\alpha$.

We assign $V$  degree 1, and $\complex G$ degree 0. This defines an
increasing filtration $F_\bullet $ on $H_{t,c}$. It was proved by in
\cite[Theorem 1.3]{eg:cherednik} that $H_{t,c}$ satisfies
Poincar\'e-Birkhoff-Witt property, i.e. the tautological embedding
$V\hookrightarrow gr(H_{t,c})$ extends to an isomorphism $Q:
\Poly(V)\rtimes G\to gr(H_{t,c})$ of vector spaces. We define $gr_i$
to be the projection from $gr(H_{t,c})$ to its $i-$th degree
component.

We define a formal deformation quantization of $\Poly(V)\rtimes
G$ as follows. For $fU_\alpha, gU_\beta\in \Poly(V)\rtimes
G$,
\[
fU_\alpha\star gU_\beta=\sum_{i,j,k=0}^\infty
\hbar^{i+j-k}Q^{-1}\Big(
gr_k\big(gr_i(Q(fU_\alpha))gr_j(Q(gU_\beta))\big)\Big).
\]
In particular, $C_i(fU_\alpha, gU_\beta)$ is defined to be
\[
C_i(fU_\alpha, gU_\beta)=\sum_{p+q-r=i}Q^{-1}\Big(
gr_r\big(gr_p(Q(fU_\alpha))gr_q(Q(gU_\beta))\big)\Big).
\]
Because $Q(fU_\alpha)$ and $Q(gU_\beta)$ are of finite degrees,
$p,q$ in the summation are both finite. Therefore, the sum in the
definition of $C_i$ is finite and $C_i$ is well defined.

We check that $\star$ is associative. $(fU_\alpha\star
gU_\beta)\star hU_{g}=$ {\small
\[
\begin{split}
=&\sum_{i,j,k=0}^\infty \hbar^{i+j-k} Q^{-1}\Big(
gr_k\big(gr_i(Q(fU_\alpha))gr_j(Q(gU_\beta))\big)\Big)\star hU_g\\
=&\sum_{i,j,k, p,q,r=0}^\infty\hbar^{i+j-k}
\hbar^{p+q-r}Q^{-1}\Big(gr_r\big(gr_p\big(Q(Q^{-1}\Big(gr_k\big(gr_i (Q(fU_\alpha)
gr_j(Q(gU_\beta)))\big)\Big))\big)\\
&\qquad \qquad gr_q(Q(hU_g))\big)\Big)\\
=&\sum_{i,j,k=0}^\infty \sum_{p,q,r=0}^\infty \hbar^{i+j+p+q-k-r}Q^{-1}
\Big(gr_r\big(gr_p\Big(gr_k\big(gr_i(Q(fU_\alpha))gr_j(Q(gU_\beta))\big)\Big)gr_q(Q(hU_g))\big)\Big)\\
=&\sum_{i,j,k=p,q,r}^\infty
\hbar^{i+j+q-r}Q^{-1}\Big(gr_r\big(gr_k\Big(gr_i(Q(fU_\alpha))gr_j(Q(gU_\beta))\Big)gr_q(Q(hU_g))\big)\Big)\\
=&\sum_{i,j,q=0}^\infty
\hbar^{i+j+q}\sum_{r=0}^\infty\hbar^{i+j+q-r}Q^{-1}\Big(gr_r\Big(\sum_{k=0}^\infty
gr_k\Big(gr_i(Q(fU_\alpha))gr_j(Q(gU_\beta))\Big)gr_q(Q(hU_g))\Big)\Big)\\
=&\sum_{i,j,k,r=0}^\infty\hbar^{i+j+q-r}Q^{-1}\Big(gr_r\Big(gr_i(Q(fU_\alpha))gr_j(Q(gU_\beta))gr_q(Q(hU_g))\Big)\Big),
\end{split}
\]}
which by the similar computation is equal to
\[
fU_\alpha\star(gU_\beta\star hU_g).
\]

We look at $C_1(fU_\alpha,
gU_\beta)=\sum_{i+j-k=1}Q^{-1}\Big(gr_k\big(gr_i(Q(fU_\alpha))gr_j(Q(gU_\beta))\big)\Big)$.
To check that $C_1$ is cohomologous to $\kappa$, we compute $L(C_1)$
as follows using the definition of $L$.
\[
\begin{split}
=&L_3\Big(\sum_{i_1, i_2}C_1\big(x_1^{i_1}-x^{i_1},
x_2^{i_2}-x^{i_2}\big)\frac{\partial}{\partial x^{i_1}}\wedge
\frac{\partial}{\partial x^{i_2}}\Big)\\
=&L_3\Big(\sum_{i_1,
i_2}gr_0\big((x_1-x)^{i_1}(x_2-x)^{i_2}\big)\frac{\partial}{\partial
x^{i_1}}\wedge
\frac{\partial}{\partial x^{i_2}}\Big)\\
=&L_3\Big(\sum_{i_2<i_1}\big(x^{i_1}x^{i_2}-x^{i_2}x^{i_1}\big)\frac{\partial}{\partial
x^{i_1}}\wedge\frac{\partial}{\partial x^{i_2}}\Big)\\
=&\sum_{i_2<i_1}\Big(\omega(x^{i_1}, x^{i_2})+\sum_{\alpha\in
S}c_\alpha\omega_\alpha(x^{i_1},
x^{i_2})U_\alpha\Big)\frac{\partial}{\partial x^{i_1}}\wedge \frac{\partial}{\partial x^{i_2}}\\
=&\frac{1}{2}\sum_{i_1,i_2}\Big(\omega(x^{i_1},
x^{i_2})+\sum_{\alpha\in S}c_\alpha\omega_\alpha(x^{i_1},
x^{i_2})U_\alpha\Big)\frac{\partial}{\partial x^{i_1}}\wedge
\frac{\partial}{\partial x^{i_2}}.
\end{split}
\]
In the third equality of the above equation, we used the
definition of $C_1$ and the product structure in $H_{t,c}$. When
$i_1<i_2$, the term $x^{i_1}x^{i_2}$ has no degree 0 term. When $i_1>i_2$,
degree 0 term of $x^{i_1}x^{i_2}$ is
$x^{i_1}x^{i_2}-x^{i_2}x^{i_1}$.

In conclusion, $\star$ is a deformation quantization of $A\rtimes
G$ with the noncommutative Poisson structure equal to
$\frac{1}{2}\kappa$.
\end{proof}

\begin{rmk}
The generalization of the deformation quantization defined in
Proposition \ref{prop:cherednik} to affine varieties was studied
by Etingof \cite{e:glob-quot}.
\end{rmk}

It is natural to ask whether the deformation quantization
constructed in Theorem \ref{prop:cherednik} is unique up to
isomorphism. From Poisson geometry, we know that
the isomorphism classes of a deformation quantization of a Poisson
structure is determined by its second Poisson cohomology, introduced in \cite{xu}.

We briefly recall the definition of Poisson cohomology here. Let $\pi$ be a noncommutative Poisson structure on an $A$. As $[\pi, \pi]_G=0$ in $HH^\bullet(A, A)$, the operator $d^\pi:HH^\bullet(A,A)\to HH^{\bullet+1}(A,A)$ defined by
\[
d^\pi(\phi)=[\pi, \phi],\qquad \qquad \phi\in H^\bullet(A,A),
\]
satisfies $d^\pi\circ d^\pi=0$. The cohomology of $HH^\bullet(A,A)$ with respect to $d^\pi$ is called the Poisson cohomology of $\pi$.

We denote the Poisson cohomology group of $Poly(V)\rtimes G$
associated to $\kappa$ by $H^\bullet _{\kappa}(Poly(V)\rtimes G)$.
In the following we compute $H^2_{\kappa}(Poly(V)\rtimes G)$.
\begin{prop}
\label{prop:poisson-coh}Let $(V, \omega)$ be a real symplectic
vector space with a finite group $G$ symplectic action. Let $\kappa$
be defined as in Corollary \ref{cor:linear}. Then the second Poisson
cohomology of $\kappa$ is isomorphic to
\[
\begin{split}
\Big\{ \sum_{g\in S}c_g\pi_g \in \Gamma^{\infty}(\bigoplus_{g\in
S}\wedge^2N^g)^G|\ \text{for all }g\in S, c_g \text{ is a constant
on $V^g$}\Big\} .
\end{split}
\]
\end{prop}
\begin{proof}
By \cite{gu-st:fixed-pt}, the fixed point subspace $V^g$ of $g$ is always a symplectic vector space, which implies that $V^g$ is of even dimension.
According to Equation (\ref{eq:cohomology}) and the fact that the
codimensions of all $V^\alpha$ are even, the second Hochschild
cohomology of $Poly(V)\rtimes G$ is
\[
\Gamma^\infty(\wedge^2 TV)^G\bigoplus (\sum_{g\in S}\wedge^0
V^g\otimes \wedge^2N^g)^G.
\]

Let $\Xi+\sum_{g\in S} f_g\pi_g$ be an element of  $ \Gamma^\infty(\wedge^2
TV)^G\bigoplus (\sum_{g\in S}\wedge^0 V^g\otimes \wedge^2 N^g)^G$.
We compute $L([T(\Xi+\sum_{g\in S} f_g\pi_g), T(\kappa)])$ as
follows. The computation is similar to computations in the proof of Theorem
\ref{thm:poisson}. There are four terms, which we deal with separately:
\[
\begin{split}
1)\ L([T(\Xi), T(\pi)]),\ \ \ \ \ \ \ \ \ &2)\ L([T(\Xi),
T(\sum_{g\in
S}c_g\pi_g)])\\
3)\ L([T(\sum_{g\in S}f_g\pi_g), T(\pi)]),\  \ \ \
\ &4)\ L([T(\sum_{g\in S}f_g\pi_g),
T(\sum_{g\in S}c_g\pi_g)]).
\end{split}
\]

The following computation follows similar reasoning to that found in
the proof of Theorem \ref{thm:poisson}.
\[
\begin{split}
1)& L([T(\Xi), T(\pi)])=[\Xi, \pi];\\
2)& L([T(\Xi), T(\sum_{g\in
S}c_g\pi_g)])=\sum_{g\in S}pr^g(c_g[\Xi,
\pi_g]|_{V^g});\\
3)& L([T(\sum_{g\in S}f_g \pi_g),
T(\pi)])=\sum_{g\in S}pr^g([f_g\pi_g, \pi]|_{V^g})=\sum_{g\in S}pr^g([f_g, \pi]\wedge \pi_g)|_{V^g};\\
4)& L([T(\sum_{g\in S}f_g\pi_g), T(\sum_{g\in
S}c_g\pi_g)])=0.
\end{split}
\]

According to the above computation,  $\Xi+\sum_{g\in
S}f_g \pi_g$ is $\kappa$ closed if and only if it
satisfies the following equations:
\begin{equation}
\label{eq:cocycle}
\begin{array}{ll}
1)& [\Xi, \pi]=0,\\
2)& pr^g(\big\{c_g[\Xi, \pi_g]+[f_g,\pi]\wedge
\pi_g \big\}|_{V^g})=0,\ \text{for all} g \in
G.
\end{array}
\end{equation}
We denote the space of solutions to the above equations by
$Z^2_\kappa$.

Next we compute the Poisson coboundary in
$(\Gamma^\infty(\wedge^2TV)\bigoplus_{\alpha\in
S}\wedge^2N^g)^G$.

According to Equation (\ref{eq:cohomology}), the cohomology $HH^1(Poly(V)\rtimes G,
Poly(V)\rtimes G)$ is a $G$-invariant vector field on $V$,
as $V^g$ for $g\ne id$ has at least codimension 2 for all $g$.
Let $X\in \Gamma^\infty(TV)^G$. We compute $L([T(\kappa), T(X)])$  using Example \ref{ex:bracket-poisson}:
\[
\begin{split}
&L([T(\kappa), T(X)])\\
=&[\pi, X]+\sum_{g\in s}pr^g(c_g [\pi_g,
X]|_{V^g}).
\end{split}
\]
We denote the space of elements in $\Gamma^\infty(\wedge^2
TV)^G\bigoplus (\sum_{g\in S}\wedge^0 V^g)^G$
of the above form by $B^2_\kappa$.

We want to find the quotient $Z^2_\kappa/B^2_\kappa$. Given
$\Xi\in \Gamma^\infty(\wedge^2TV)^G$ with $[\Xi, \pi]=0$,
we use the fact that $\pi$ is from a symplectic form to find a
$G$-invariant vector field $X\in \Gamma^\infty(TV)$ such that
$[\pi, X]=\Xi$. Because $L([T(\kappa), T(X)])=:\Xi+\sum_g h_g\pi_g=[\pi, X]+\sum_{g\in
S}c_g pr^g([ \pi_g, X]|_{V^g})$ is
$\kappa$-closed we conclude that
\[
\Xi+\sum_{g\in S}h_g\pi_g-(\Xi+\sum_g f_g\pi_g)=\sum_{g\in S} (h_g-f_g)\pi_g
\]
is also $\kappa$ closed. Substituting this expression into second
equation of (\ref{eq:cocycle}) with $\Xi=0$, we have that for all $g\in S$
\[
pr^g([\pi,
(h_g-f_g)]\wedge\pi_g|_{V^g})=[\pi,
(h_g-f_g)]|_{V^g}\wedge \pi_g|_{V^g}=0.
\]
In the above equality, we have used the fact that $h_g-f_g$ is supported on $V^g$ and therefore $[\pi, h_g-f_g]$ is an element in $TV^g$ also  supported on $V^g$. Accordingly, we have
\[
[\pi, h_g-f_g]|_{V^g}=0\Leftrightarrow[\pi,
h_g-f_g]=0.
\]
Because $\pi$ is symplectic, this implies that $h_g-f_g$
has to be a constant. It is obvious that $X+\sum_{g\in
S}(h_g+c_g)\pi_g$ is $\kappa$-closed.

We are left to show that nonzero elements of the form $\sum_{g\in
S}a_g\pi_g$ are not coboundaries for any constant
$a_g$, $g\in S$. If $X\in \Gamma^\infty(TV)^G$ such
that $L(T(\kappa), T(X))=\sum_{g\in S}a_g\pi_g$,
then by similar computations to those above, we have that
\[
\begin{split}
[\pi,X]=0,\ \ \ \ \ \ \ & pr^g(c_g[\pi_g,
X]|_{V^g})=a_g\pi_g.
\end{split}
\]
for any $g\in S$.

As $\pi$ is from a symplectic form, $[\pi, X]=0$ implies that there
is a function $f$ such that $X=[\pi, f]$. As $g$ acts on $V$
preserving the symplectic form, we are allowed to write $\pi$  as a
sum of $\pi_g$ and $\pi^g$ where $\pi_g$ is a bivector supported along $N^g$
and $\pi^g$ which is a bivector supported along $V^g$. Therefore, we
have
\[
[\pi_g, X]=[\pi_g, [\pi, f]]=[\pi_g, [\pi_g,f]]+[\pi_g, [\pi^g, f]].
\]
We easily compute that $[\pi_g, [\pi_g, f]]\equiv 0$ and
$pr^g([\pi_g, [\pi^g, f]])$ using the fact $N^g$ is of dimension 2.
This shows that $[\pi_g, X]$ contains no component proportional to
$\pi_g$. Accordingly, $pr^g(c_g[\pi_g, X])=0=a_g\pi_g. $

In conclusion, we see that the quotient space
$Z^2_\kappa/B^2_\kappa$ is
\[
\{ \sum_{g\in S}c_g\pi_g \in
\Gamma^{\infty}(\bigoplus_{g\in S}\wedge^2N^g)^G|\
\text{for all }g\in S, c_g \text{ is a constant on
$V^g$}.\}
\]
\end{proof}

\begin{rmk}
Proposition \ref{prop:poisson-coh} shows that the dimension of all infinitesimal deformations of a noncommutative Poisson
structure $\kappa$ is equal to the size of the set $S$.
Furthermore, it is easy to check that the all infinitesimal
deformations actually correspond to Poisson structures. This gives
another explanation of Corollary \ref{cor:linear}.
\end{rmk}

Propositions \ref{prop:cherednik} and \ref{prop:poisson-coh} inspire a
series of interesting questions. The cocycle $\kappa$ is a very
special type of noncommutative Poisson structure on $A\rtimes G$
defined in Theorem \ref{thm:poisson}. In \cite{hot:example}, with
Jean-Michel Oudom, we proved that with a mild assumption all linear
Poisson structures on $Poly(V)\rtimes G$ can be quantized. This
inspires the question of whether all the noncommutative Poisson
structures defined in Theorem \ref{thm:poisson} can be deformation
quantized. If their deformation quantizations exist, how many are
there? All these problems have a general version on
$C^\infty(M)\rtimes G$. It is closely related to the Conjecture 1 in
\cite{do-et:cohomology} by Dolgushev and Etingof. We plan to address
these questions in the near future.

\subsection{Noncommutative quadratic Poisson structures}
In this subsection, we provide some new examples of noncommutative
Poisson structures  other than those in Corollary
\ref{cor:linear}. We can easily see that these Poisson structures
are not symplectic at all, and they should be viewed as
generalized quadratic Poisson structures.

We consider the space of $\reals^4=\complex \times \complex$ with
the following $\integers_n\times \integers_m$ action, where
$\integers_n=\integers/n\integers$ and
$\integers_m=\integers/m\integers$. Let $(z_1, z_2)$ be
holomorphic coordinates  on $\complex^2$, and $(k,l)\in
\integers_n\times \integers_m$. Define
\[
\begin{split}
(k,l):&(z_1, z_2)\longrightarrow (\exp(\frac{2k\pi i}{n})z_1, \exp(\frac{2l \pi i}{m})z_2)\\
&(\bar{z}_1, \bar{z}_2)\longrightarrow (\exp(-\frac{2k\pi
i}{n})\bar{z}_1, \exp(-\frac{2l\pi i}{m})\bar{z}_2).
\end{split}
\]

The fixed point subspace of $(k,l)\in \integers_n\times
\integers_m$ can be described explicitly.
\begin{enumerate}
\item if $k\ne 0, l\ne0$, the fixed point set of $(k,l)$
consists of only one point, the origin;
\item if $k=0, l\ne0$, $(0,l)$'s fixed point set is
$\complex\times \{0\}\subset \complex \times \complex$;
\item if $k\ne 0, l=0$, $(k,0)$'s fixed point set is
$\{0\}\times \complex\subset \complex \times \complex$;
\item if $k=l=0$, $(0,0)$'s fixed point space is
$\complex \times \complex$.
\end{enumerate}

To look for noncommutative Poisson structures on
$Poly(\reals^4)\rtimes (\integers_n\times \integers_m)$, we only
need to consider the fixed point space of the identity, which is
$\complex^2$, and the codimension 2 fixed point subspace of the form $\{0\}\times \complex$ of $(k,0)$, and
the fixed point subspace $\complex\times \{0\}$ of $(0,l)$. We consider
the following collection of bivector fields where $\alpha,\beta,
\lambda_k, \mu_l$ are real constants.
\begin{enumerate}
\item on the fixed point subspace of $(0,0)$, we consider
$\Pi_{0,0}^\alpha=i\alpha |z_2|^2\frac{\partial}{\partial
z_1}\wedge \frac{\partial}{\partial\bar{z}_1}$. We notice that
$\Pi_{0,0}$ is $\integers_n\times \integers_m$ invariant, and
satisfies $[\Pi_{0,0}, \Pi_{0,0}]=0$;
\item on  the fixed point subspace of $(k,0)$, we consider $\Pi_{k,0}=
i\lambda_k |z_2|^2\frac{\partial}{\partial z_1}\wedge
\frac{\partial}{\partial \bar{z}_1}$, which is a smooth section of
the determinant bundle of the normal bundle over $\{0\}\times
\complex\subset \complex\times \complex$. Again we notice that
$\Pi_{k,0}$ is $\integers_n\times\integers_m$ invariant, and
$[\Pi_{0,0}, \Pi_{k,0}]=0$;
\item on the fixed point subspace of $(0,l)$, we consider $\Pi_{0,l}=i\mu_l
|z_1|^2\frac{\partial}{\partial z_2}\wedge
\frac{\partial}{\partial\bar{z}_2}$, which is a smooth section of
the determinant bundle of the normal bundle over $\complex \times
\{0\}\subset \complex\times \complex$. We notice that $\Pi_{0,l}$
is $\integers_n\times \integers_m$ invariant, but
\[
\begin{split}
[\Pi_{0,0},
\Pi_{0,l}]&=-\alpha\mu_l\Big[|z_2|^2(z_1\frac{\partial}{\partial
z_1} -\bar{z}_1\frac{\partial}{\partial\bar{z}_1})\wedge
\frac{\partial}{\partial z_2}\wedge \frac{\partial}{\partial\bar{z}_2}\\
&+|z_1|^2(z_2\frac{\partial}{\partial
z_2}-\bar{z}_2\frac{\partial}{\partial\bar{z}_2})\wedge
\frac{\partial}{\partial z_1}\wedge
\frac{\partial}{\partial\bar{z}_1}\Big]\ne 0.
\end{split}
\]

However, if we look at the restriction of $[\Pi_{0,0}, \Pi_{0,l}]$
to the fixed point subspace of $(0,l)$ which is $\complex\times \{0\}$,
it does vanish. We have that
\[
pr^{(0,l)}([\Pi_{0,0}, \Pi_{0,l}]|_{\complex\times \{0\}})=0.
\]
\end{enumerate}

By Theorem \ref{thm:poisson}, we conclude that the collection of
$(\Pi^\alpha_{0,0}, \Pi_{k,0},\Pi_{0,l})$ defines a family of
noncommutative Poisson structures on $Poly(\reals^4)\rtimes
(\integers_n\times \integers_m)$. One can also easily check that
if we replace $\Pi^\alpha_{0,0}$ by
$\Pi^\beta_{0,0}=i\beta|z_1|^2\frac{\partial}{\partial z_2}\wedge
\frac{\partial}{\partial\bar{z}_2}$, ($\Pi_{0,0}^\beta, \Pi_{k,0},
\Pi_{0,l})$ defines another family of noncommutative Poisson
structures on $Poly(\reals^4)\rtimes (\integers_n\times
\integers_m)$.

There are various ways to generalize the above families of
examples to higher dimensions. For example, one can consider actions of
$\integers_n\times \integers_m$ which act on the first two components
of $\complex^k$ as above, but act trivially on the
left $\complex^{k-2}$ component. Then the above two families of
noncommutative Poisson structures naturally extend to
$Poly(\reals^{2k})\rtimes(\integers_n\times \integers_m)$. We will
leave the more nontrivial generalizations to the future
\cite{hot:example}.

\bibliographystyle{amsplain}

\end{document}